\documentclass[11pt,a4paper]{article}

\usepackage{amsmath,amssymb,latexsym,pstricks,mathrsfs,comment,amsthm,graphicx,tikz,tikz-cd,enumerate,accents,pgffor,cite,wrapfig,multicol,float,slashed}
\usepackage{bibspacing}
\usepackage[T1]{fontenc}
\usepackage{ wasysym }
\usepackage{stmaryrd}
\usepackage[colorlinks]{hyperref}

\newcommand{\nc}{\newcommand}
\nc{\rnc}{\renewcommand}

\usepackage{geometry} 

\let\oldproofname=\proofname
\renewcommand{\proofname}{\rm\bf{\oldproofname}}

\makeatletter

\DeclareMathSymbol{\widehatsym}{\mathord}{largesymbols}{"62}
\newcommand\lowerwidehatsym{%
  \text{\smash{\raisebox{-1.3ex}{%
    $\widehatsym$}}}}
\newcommand\fixwidehat[1]{%
  \mathchoice
    {\accentset{\displaystyle\lowerwidehatsym}{#1}}
    {\accentset{\textstyle\lowerwidehatsym}{#1}}
    {\accentset{\scriptstyle\lowerwidehatsym}{#1}}
    {\accentset{\scriptscriptstyle\lowerwidehatsym}{#1}}
}
\rnc{\widehat}{\fixwidehat}

\begin{document}

\nc{\pfitem}[1]{\medskip \noindent (#1).}
\nc{\Set}{\operatorname{{\bf Set}}}
\nc{\Setp}{\operatorname{{\bf Set}}^+}
\nc{\TXYa}{\T_{XY}^a}
\nc{\IXYa}{\I_{XY}^a}
\nc{\PT}{\P\T}
\rnc{\emptyset}{\varnothing}
\nc{\leqJa}{\leq_{\J^a}}
\nc{\leqRa}{\leq_{\R^a}}
\nc{\leqLa}{\leq_{\L^a}}
\rnc{\th}{\theta}
\nc{\starb}{\oast}
\nc{\bPTfr}{\bPT^{\operatorname{fr}}}
\nc{\Sym}{\operatorname{Sym}}
\nc{\PTXYa}{\PT_{XY}^a}
\nc{\PTXY}{\PT_{XY}}
\nc{\TXY}{\T_{XY}}
\nc{\IXY}{\I_{XY}}
\nc{\PTYX}{\PT_{YX}}
\nc{\tran}[2]{\left(\begin{smallmatrix} #1\\#2 \end{smallmatrix}\right)}
\nc{\gt}{\widetilde{g}}
\nc{\cE}{\mathcal E}
\nc{\EXYa}{\cE_{XY}^a}
\nc{\stirlingii}{\genfrac{[}{]}{0pt}{}}
\nc{\sh}{\operatorname{sh}}
\nc{\col}{\operatorname{col}}
\nc{\defect}{\operatorname{def}}
\nc{\codefect}{\operatorname{codef}}
\nc{\sP}{\mathscr P}
\nc{\sE}{\mathscr E}
\nc{\sF}{\mathscr F}
\nc{\JPa}{\J^{P^a}}
\nc{\leqJPa}{\leq_{\JPa}}
\nc{\A}{B} 
\nc{\B}{A} 
\nc{\sectiontitle}[1]{\section{\boldmath #1}}
\nc{\subsectiontitle}[1]{\subsection[#1]{\boldmath #1}}

\nc{\MaxE}{\operatorname{Max}_\preceq}

\nc{\DClass}[5]{
\begin{tikzpicture}[scale=#5]
\foreach \x/\y in {#3} {\fill[lightgray](\y-1,#2-\x)--(\y,#2-\x)--(\y,#2-\x+1)--(\y-1,#2-\x+1)--(\y-1,#2-\x); \draw(\y-.5,#2-\x+.5)node{$\scriptscriptstyle{#4}$};}
\foreach \x in {0,1,...,#1} {\draw(\x,0)--(\x,#2);}
\foreach \x in {0,1,...,#2} {\draw(0,\x)--(#1,\x);}
\draw[ultra thick] (0,0)--(#1,0)--(#1,#2)--(0,#2)--(0,0)--(#1,0);
\end{tikzpicture}
}

\nc{\DaClass}[7]{
\begin{tikzpicture}[scale=#5]
\foreach \x/\y in {#3} {\fill[lightgray](\y-1,#2-\x)--(\y,#2-\x)--(\y,#2-\x+1)--(\y-1,#2-\x+1)--(\y-1,#2-\x); \draw(\y-.5,#2-\x+.5)node{$\scriptscriptstyle{#4}$};}
\foreach \x in {0,1,...,#1} {\draw(\x,0)--(\x,#2);}
\foreach \x in {0,1,...,#2} {\draw(0,\x)--(#1,\x);}
\foreach \x in {#6} {\draw[ultra thick](\x,0)--(\x,#2);}
\foreach \x in {#7} {\draw[ultra thick](0,\x)--(#1,\x);}
\draw[ultra thick] (0,0)--(#1,0)--(#1,#2)--(0,#2)--(0,0)--(#1,0);
\end{tikzpicture}
}

\nc{\bT}{\operatorname{\bf T}}
\nc{\bPT}{\operatorname{\bf PT}}
\nc{\bI}{\operatorname{\bf I}}

\nc{\M}{\mathcal M}
\nc{\G}{\mathcal G}
\nc{\F}{\mathbb F}
\nc{\MnJ}{\mathcal M_n^J}
\nc{\EnJ}{\mathcal E_n^J}
\nc{\Mat}{\operatorname{Mat}}
\nc{\RegMnJ}{\Reg(\MnJ)}
\nc{\row}{\mathfrak r}
\nc{\Row}{\operatorname{Row}}
\nc{\Col}{\operatorname{Col}}
\nc{\Span}{\operatorname{span}}
\nc{\mat}[4]{\left[\begin{matrix}#1&#2\\#3&#4\end{matrix}\right]}
\nc{\tmattwo}[2]{\left[\begin{smallmatrix}#1\\#2\end{smallmatrix}\right]}
\nc{\tmat}[4]{\left[\begin{smallmatrix}#1&#2\\#3&#4\end{smallmatrix}\right]}
\nc{\ttmat}[4]{{\tiny \left[\begin{smallmatrix}#1&#2\\#3&#4\end{smallmatrix}\right]}}
\nc{\tmatt}[9]{\left[\begin{smallmatrix}#1&#2&#3\\#4&#5&#6\\#7&#8&#9\end{smallmatrix}\right]}
\nc{\ttmatt}[9]{{\tiny \left[\begin{smallmatrix}#1&#2&#3\\#4&#5&#6\\#7&#8&#9\end{smallmatrix}\right]}}
\nc{\MnGn}{\M_n\sm\G_n}
\nc{\MrGr}{\M_r\sm\G_r}
\nc{\qbin}[2]{\left[\begin{matrix}#1\\#2\end{matrix}\right]_q}
\nc{\tqbin}[2]{\left[\begin{smallmatrix}#1\\#2\end{smallmatrix}\right]_q}
\nc{\qbinx}[3]{\left[\begin{matrix}#1\\#2\end{matrix}\right]_{#3}}
\nc{\tqbinx}[3]{\left[\begin{smallmatrix}#1\\#2\end{smallmatrix}\right]_{#3}}
\nc{\MNJ}{\M_nJ}
\nc{\JMN}{J\M_n}
\nc{\RegMNJ}{\Reg(\MNJ)}
\nc{\RegJMN}{\Reg(\JMN)}
\nc{\RegMMNJ}{\Reg(\MMNJ)}
\nc{\RegJMMN}{\Reg(\JMMN)}
\nc{\Wb}{\overline{W}}
\nc{\Xb}{\overline{X}}
\nc{\Yb}{\overline{Y}}
\nc{\Zb}{\overline{Z}}
\nc{\Sib}{\overline{\Si}}
\nc{\Om}{\Omega}
\nc{\Omb}{\overline{\Om}}
\nc{\Gab}{\overline{\Ga}}
\nc{\qfact}[1]{[#1]_q!}
\nc{\smat}[2]{\left[\begin{matrix}#1&#2\end{matrix}\right]}
\nc{\tsmat}[2]{\left[\begin{smallmatrix}#1&#2\end{smallmatrix}\right]}
\nc{\hmat}[2]{\left[\begin{matrix}#1\\#2\end{matrix}\right]}
\nc{\thmat}[2]{\left[\begin{smallmatrix}#1\\#2\end{smallmatrix}\right]}
\nc{\LVW}{\mathcal L(V,W)}
\nc{\KVW}{\mathcal K(V,W)}
\nc{\LV}{\mathcal L(V)}
\nc{\RegLVW}{\Reg(\LVW)}
\nc{\sM}{\mathscr M}
\nc{\sN}{\mathscr N}
\rnc{\iff}{\ \Leftrightarrow\ }
\nc{\Hom}{\operatorname{Hom}}
\nc{\End}{\operatorname{End}}
\nc{\Aut}{\operatorname{Aut}}
\nc{\Lin}{\mathcal L}
\nc{\Hommn}{\Hom(V_m,V_n)}
\nc{\Homnm}{\Hom(V_n,V_m)}
\nc{\Homnl}{\Hom(V_n,V_l)}
\nc{\Homkm}{\Hom(V_k,V_m)}
\nc{\Endm}{\End(V_m)}
\nc{\Endn}{\End(V_n)}
\nc{\Endr}{\End(V_r)}
\nc{\Autm}{\Aut(V_m)}
\nc{\Autn}{\Aut(V_n)}
\nc{\MmnJ}{\M_{mn}^J}
\nc{\MmnA}{\M_{mn}^A}
\nc{\MmnB}{\M_{mn}^B}
\nc{\Mmn}{\M_{mn}}
\nc{\Mkl}{\M_{kl}}
\nc{\Mnm}{\M_{nm}}
\nc{\EmnJ}{\mathcal E_{mn}^J}
\nc{\MmGm}{\M_m\sm\G_m}
\nc{\RegMmnJ}{\Reg(\MmnJ)}
\rnc{\implies}{\ \Rightarrow\ }
\nc{\DMmn}[1]{D_{#1}(\Mmn)}
\nc{\DMmnJ}[1]{D_{#1}(\MmnJ)}
\nc{\MMNJ}{\Mmn J}
\nc{\JMMN}{J\Mmn}
\nc{\JMMNJ}{J\Mmn J}
\nc{\Inr}{\mathcal I(V_n,W_r)}
\nc{\Lnr}{\mathcal L(V_n,W_r)}
\nc{\Knr}{\mathcal K(V_n,W_r)}
\nc{\Imr}{\mathcal I(V_m,W_r)}
\nc{\Kmr}{\mathcal K(V_m,W_r)}
\nc{\Lmr}{\mathcal L(V_m,W_r)}
\nc{\Kmmr}{\mathcal K(V_m,W_{m-r})}
\nc{\tr}{{\operatorname{T}}}
\nc{\MMN}{\MmnA(\F_1)}
\nc{\MKL}{\Mkl^B(\F_2)}
\nc{\RegMMN}{\Reg(\MmnA(\F_1))}
\nc{\RegMKL}{\Reg(\Mkl^B(\F_2))}
\nc{\gRhA}{\widehat{\mathscr R}^A}
\nc{\gRhB}{\widehat{\mathscr R}^B}
\nc{\gLhA}{\widehat{\mathscr L}^A}
\nc{\gLhB}{\widehat{\mathscr L}^B}
\nc{\timplies}{\Rightarrow}
\nc{\tiff}{\Leftrightarrow}
\nc{\Sija}{S_{ij}^a}
\nc{\dmat}[8]{\draw(#1*1.5,#2)node{$\left[\begin{smallmatrix}#3&#4&#5\\#6&#7&#8\end{smallmatrix}\right]$};}
\nc{\bdmat}[8]{\draw(#1*1.5,#2)node{${\mathbf{\left[\begin{smallmatrix}#3&#4&#5\\#6&#7&#8\end{smallmatrix}\right]}}$};}
\nc{\rdmat}[8]{\draw(#1*1.5,#2)node{\rotatebox{90}{$\left[\begin{smallmatrix}#3&#4&#5\\#6&#7&#8\end{smallmatrix}\right]$}};}
\nc{\rldmat}[8]{\draw(#1*1.5-0.375,#2)node{\rotatebox{90}{$\left[\begin{smallmatrix}#3&#4&#5\\#6&#7&#8\end{smallmatrix}\right]$}};}
\nc{\rrdmat}[8]{\draw(#1*1.5+.375,#2)node{\rotatebox{90}{$\left[\begin{smallmatrix}#3&#4&#5\\#6&#7&#8\end{smallmatrix}\right]$}};}
\nc{\rfldmat}[8]{\draw(#1*1.5-0.375+.15,#2)node{\rotatebox{90}{$\left[\begin{smallmatrix}#3&#4&#5\\#6&#7&#8\end{smallmatrix}\right]$}};}
\nc{\rfrdmat}[8]{\draw(#1*1.5+.375-.15,#2)node{\rotatebox{90}{$\left[\begin{smallmatrix}#3&#4&#5\\#6&#7&#8\end{smallmatrix}\right]$}};}
\nc{\xL}{[x]_{\! _\gL}}\nc{\yL}{[y]_{\! _\gL}}\nc{\xR}{[x]_{\! _\gR}}\nc{\yR}{[y]_{\! _\gR}}\nc{\xH}{[x]_{\! _\gH}}\nc{\yH}{[y]_{\! _\gH}}\nc{\XK}{[X]_{\! _\gK}}\nc{\xK}{[x]_{\! _\gK}}
\nc{\RegSija}{\Reg(\Sija)}
\nc{\MnmK}{\M_{nm}^K}
\nc{\cC}{\mathcal C}
\nc{\cR}{\mathcal R}
\nc{\Ckl}{\cC_k(l)}
\nc{\Rkl}{\cR_k(l)}
\nc{\Cmr}{\cC_m(r)}
\nc{\Rmr}{\cR_m(r)}
\nc{\Cnr}{\cC_n(r)}
\nc{\Rnr}{\cR_n(r)}
\nc{\Z}{\mathbb Z}

\nc{\Reg}{\operatorname{Reg}}
\nc{\RP}{\operatorname{RP}}
\nc{\MI}{\operatorname{MI}}
\nc{\TXa}{\T_X^a}
\nc{\TXA}{\T(X,A)}
\nc{\TXal}{\T(X,\al)}
\nc{\RegTXa}{\Reg(\TXa)}
\nc{\RegTXA}{\Reg(\TXA)}
\nc{\RegTXal}{\Reg(\TXal)}
\nc{\PalX}{\P_\al(X)}
\nc{\EAX}{\E_A(X)}
\nc{\Bb}{\overline{B}}
\nc{\bb}{\overline{b}}
\nc{\bw}{{\bf w}}
\nc{\bz}{{\bf z}}
\nc{\TASA}{\T_A\sm\S_A}
\nc{\Ub}{\overline{U}}
\nc{\Vb}{\overline{V}}
\nc{\eb}{\overline{e}}
\nc{\EXa}{\E_X^a}
\nc{\oijr}{1\leq i<j\leq r}
\nc{\veb}{\overline{\ve}}
\nc{\bbT}{\mathbb T}
\nc{\Surj}{\operatorname{Surj}}
\nc{\Sone}{S^{(1)}}
\nc{\fillbox}[2]{\draw[fill=gray!30](#1,#2)--(#1+1,#2)--(#1+1,#2+1)--(#1,#2+1)--(#1,#2);}
\nc{\raa}{\rangle_J}
\nc{\raJ}{\rangle_J}
\nc{\Ea}{E_J}
\nc{\EJ}{E_J}
\nc{\ep}{\epsilon} \nc{\ve}{\varepsilon}
\nc{\IXa}{\I_X^a}
\nc{\RegIXa}{\Reg(\IXa)}
\nc{\JXa}{\J_X^a}
\nc{\RegJXa}{\Reg(\JXa)}
\nc{\IXA}{\I(X,A)}
\nc{\IAX}{\I(A,X)}
\nc{\RegIXA}{\Reg(\IXA)}
\nc{\RegIAX}{\Reg(\IAX)}
\nc{\trans}[2]{\left(\begin{smallmatrix} #1 \\ #2 \end{smallmatrix}\right)}
\nc{\bigtrans}[2]{\left(\begin{matrix} #1 \\ #2 \end{matrix}\right)}
\nc{\lmap}[1]{\mapstochar \xrightarrow {\ #1\ }}
\nc{\EaTXa}{E}

\nc{\gL}{\mathrel{\mathscr L}}
\nc{\gR}{\mathrel{\mathscr R}}
\nc{\gH}{\mathrel{\mathscr H}}
\nc{\gJ}{\mathrel{\mathscr J}}
\nc{\gD}{\mathrel{\mathscr D}}
\nc{\gK}{\mathrel{\mathscr K}}
\nc{\gLa}{\mathrel{\mathscr L^a}}
\nc{\gRa}{\mathrel{\mathscr R^a}}
\nc{\gHa}{\mathrel{\mathscr H^a}}
\nc{\gJa}{\mathrel{\mathscr J^a}}
\nc{\gDa}{\mathrel{\mathscr D^a}}
\nc{\gKa}{\mathrel{\mathscr K^a}}
\nc{\gLJ}{\mathrel{\mathscr L^J}}
\nc{\gRJ}{\mathrel{\mathscr R^J}}
\nc{\gHJ}{\mathrel{\mathscr H^J}}
\nc{\gJJ}{\mathrel{\mathscr J^J}}
\nc{\gDJ}{\mathrel{\mathscr D^J}}
\nc{\gKJ}{\mathrel{\mathscr K^J}}
\nc{\gLh}{\mathrel{\widehat{\mathscr L}}}
\nc{\gRh}{\mathrel{\widehat{\mathscr R}}}
\nc{\gHh}{\mathrel{\widehat{\mathscr H}}}
\nc{\gJh}{\mathrel{\widehat{\mathscr J}}}
\nc{\gDh}{\mathrel{\widehat{\mathscr D}}}
\nc{\gKh}{\mathrel{\widehat{\mathscr K}}}
\nc{\Lh}{\widehat{L}}
\nc{\Rh}{\widehat{R}}
\nc{\Hh}{\widehat{H}}
\nc{\Jh}{\widehat{J}}
\nc{\Dh}{\widehat{D}}
\nc{\Kh}{\widehat{K}}
\nc{\gLb}{\mathrel{\widehat{\mathscr L}}}
\nc{\gRb}{\mathrel{\widehat{\mathscr R}}}
\nc{\gHb}{\mathrel{\widehat{\mathscr H}}}
\nc{\gJb}{\mathrel{\widehat{\mathscr J}}}
\nc{\gDb}{\mathrel{\widehat{\mathscr D}}}
\nc{\gKb}{\mathrel{\widehat{\mathscr K}}}
\nc{\Lb}{\mathrel{\widehat{L}^J}}
\nc{\Rb}{\mathrel{\widehat{R}^J}}
\nc{\Hb}{\mathrel{\widehat{H}^J}}
\nc{\Jb}{\mathrel{\widehat{J}^J}}
\nc{\Db}{\mathrel{\overline{D}}}
\nc{\Kb}{\mathrel{\widehat{K}}}

\nc{\xb}{\overline{x}}
\nc{\yb}{\overline{y}}
\nc{\zb}{\overline{z}}
\nc{\qb}{\overline{q}}

\hyphenation{mon-oid mon-oids}

\nc{\itemit}[1]{\item[\emph{(#1)}]}
\nc{\itemnit}[1]{\item[(#1)]}
\nc{\E}{\mathbb E}
\nc{\TX}{\T(X)}
\nc{\TXP}{\T(X,\P)}
\nc{\EX}{\E(X)}
\nc{\EXP}{\E(X,\P)}
\nc{\SX}{\S(X)}
\nc{\SXP}{\S(X,\P)}
\nc{\Sing}{\operatorname{Sing}}
\nc{\idrank}{\operatorname{idrank}}
\nc{\SingXP}{\Sing(X,\P)}
\nc{\De}{\Delta}
\nc{\sgp}{\operatorname{sgp}}
\nc{\mon}{\operatorname{mon}}
\nc{\Dn}{\mathcal D_n}
\nc{\Dm}{\mathcal D_m}

\nc{\lline}[1]{\draw(3*#1,0)--(3*#1+2,0);}
\nc{\uline}[1]{\draw(3*#1,5)--(3*#1+2,5);}
\nc{\thickline}[2]{\draw(3*#1,5)--(3*#2,0); \draw(3*#1+2,5)--(3*#2+2,0) ;}
\nc{\thicklabel}[3]{\draw(3*#1+1+3*#2*0.15-3*#1*0.15,4.25)node{{\tiny $#3$}};}

\nc{\slline}[3]{\draw(3*#1+#3,0+#2)--(3*#1+2+#3,0+#2);}
\nc{\suline}[3]{\draw(3*#1+#3,5+#2)--(3*#1+2+#3,5+#2);}
\nc{\sthickline}[4]{\draw(3*#1+#4,5+#3)--(3*#2+#4,0+#3); \draw(3*#1+2+#4,5+#3)--(3*#2+2+#4,0+#3) ;}
\nc{\sthicklabel}[5]{\draw(3*#1+1+3*#2*0.15-3*#1*0.15+#5,4.25+#4)node{{\tiny $#3$}};}

\nc{\stll}[5]{\sthickline{#1}{#2}{#4}{#5} \sthicklabel{#1}{#2}{#3}{#4}{#5}}
\nc{\tll}[3]{\stll{#1}{#2}{#3}00}

\nc{\mfourpic}[9]{
\slline1{#9}0
\slline3{#9}0
\slline4{#9}0
\slline5{#9}0
\suline1{#9}0
\suline3{#9}0
\suline4{#9}0
\suline5{#9}0
\stll1{#1}{#5}{#9}{0}
\stll3{#2}{#6}{#9}{0}
\stll4{#3}{#7}{#9}{0}
\stll5{#4}{#8}{#9}{0}
\draw[dotted](6,0+#9)--(8,0+#9);
\draw[dotted](6,5+#9)--(8,5+#9);
}
\nc{\vdotted}[1]{
\draw[dotted](3*#1,10)--(3*#1,15);
\draw[dotted](3*#1+2,10)--(3*#1+2,15);
}

\nc{\Clab}[2]{
\sthicklabel{#1}{#1}{{}_{\phantom{#1}}C_{#1}}{1.25+5*#2}0
}
\nc{\sClab}[3]{
\sthicklabel{#1}{#1}{{}_{\phantom{#1}}C_{#1}}{1.25+5*#2}{#3}
}
\nc{\Clabl}[3]{
\sthicklabel{#1}{#1}{{}_{\phantom{#3}}C_{#3}}{1.25+5*#2}0
}
\nc{\sClabl}[4]{
\sthicklabel{#1}{#1}{{}_{\phantom{#4}}C_{#4}}{1.25+5*#2}{#3}
}
\nc{\Clabll}[3]{
\sthicklabel{#1}{#1}{C_{#3}}{1.25+5*#2}0
}
\nc{\sClabll}[4]{
\sthicklabel{#1}{#1}{C_{#3}}{1.25+5*#2}{#3}
}

\nc{\mtwopic}[6]{
\slline1{#6*5}{#5}
\slline2{#6*5}{#5}
\suline1{#6*5}{#5}
\suline2{#6*5}{#5}
\stll1{#1}{#3}{#6*5}{#5}
\stll2{#2}{#4}{#6*5}{#5}
}
\nc{\mtwopicl}[6]{
\slline1{#6*5}{#5}
\slline2{#6*5}{#5}
\suline1{#6*5}{#5}
\suline2{#6*5}{#5}
\stll1{#1}{#3}{#6*5}{#5}
\stll2{#2}{#4}{#6*5}{#5}
\sClabl1{#6}{#5}{i}
\sClabl2{#6}{#5}{j}
}

\nc{\keru}{\operatorname{ker}^\wedge} \nc{\kerl}{\operatorname{ker}_\vee}

\nc{\coker}{\operatorname{coker}}
\nc{\KER}{\ker}
\nc{\N}{\mathbb N}
\nc{\LaBn}{L_\al(\B_n)}
\nc{\RaBn}{R_\al(\B_n)}
\nc{\LaPBn}{L_\al(\PB_n)}
\nc{\RaPBn}{R_\al(\PB_n)}
\nc{\rhorBn}{\rho_r(\B_n)}
\nc{\DrBn}{D_r(\B_n)}
\nc{\DrPn}{D_r(\P_n)}
\nc{\DrPBn}{D_r(\PB_n)}
\nc{\DrKn}{D_r(\K_n)}
\nc{\alb}{\al_{\vee}}
\nc{\beb}{\be^{\wedge}}
\nc{\bnf}{\bn^\flat}
\nc{\Bal}{\operatorname{Bal}}
\nc{\Red}{\operatorname{Red}}
\nc{\Pnxi}{\P_n^\xi}
\nc{\Bnxi}{\B_n^\xi}
\nc{\PBnxi}{\PB_n^\xi}
\nc{\Knxi}{\K_n^\xi}
\nc{\C}{\mathscr C}
\nc{\exi}{e^\xi}
\nc{\Exi}{E^\xi}
\nc{\eximu}{e^\xi_\mu}
\nc{\Eximu}{E^\xi_\mu}
\nc{\REF}{ {\red [Ref?]} }
\nc{\GL}{\operatorname{GL}}
\rnc{\O}{\operatorname{O}}

\nc{\vtx}[2]{\fill (#1,#2)circle(.2);}
\nc{\lvtx}[2]{\fill (#1,0)circle(.2);}
\nc{\uvtx}[2]{\fill (#1,1.5)circle(.2);}

\nc{\Eq}{\mathfrak{Eq}}
\nc{\Gau}{\Ga^\wedge} \nc{\Gal}{\Ga_\vee}
\nc{\Lamu}{\Lam^\wedge} \nc{\Laml}{\Lam_\vee}
\nc{\bX}{{\bf X}}
\nc{\bY}{{\bf Y}}
\nc{\ds}{\displaystyle}

\nc{\uvert}[1]{\fill (#1,1.5)circle(.2);}
\nc{\uuvert}[1]{\fill (#1,3)circle(.2);}
\nc{\uuuvert}[1]{\fill (#1,4.5)circle(.2);}
\rnc{\lvert}[1]{\fill (#1,0)circle(.2);}
\nc{\overtl}[3]{\node[vertex] (#3) at (#1,0) {  {\tiny $#2$} };}
\nc{\cv}[2]{\draw(#1,1.5) to [out=270,in=90] (#2,0);}
\nc{\cvs}[2]{\draw(#1,1.5) to [out=270+30,in=90+30] (#2,0);}
\nc{\ucv}[2]{\draw(#1,3) to [out=270,in=90] (#2,1.5);}
\nc{\uucv}[2]{\draw(#1,4.5) to [out=270,in=90] (#2,3);}
\nc{\textpartn}[1]{{\lower0.45 ex\hbox{\begin{tikzpicture}[xscale=.2,yscale=0.2] #1 \end{tikzpicture}}}}
\nc{\textpartnx}[2]{{\lower1.0 ex\hbox{\begin{tikzpicture}[xscale=.3,yscale=0.3] 
\foreach \x in {1,...,#1}
{ \uvert{\x} \lvert{\x} }
#2 \end{tikzpicture}}}}
\nc{\disppartnx}[2]{{\lower1.0 ex\hbox{\begin{tikzpicture}[scale=0.3] 
\foreach \x in {1,...,#1}
{ \uvert{\x} \lvert{\x} }
#2 \end{tikzpicture}}}}
\nc{\disppartnxd}[2]{{\lower2.1 ex\hbox{\begin{tikzpicture}[scale=0.3] 
\foreach \x in {1,...,#1}
{ \uuvert{\x} \uvert{\x} \lvert{\x} }
#2 \end{tikzpicture}}}}
\nc{\disppartnxdn}[2]{{\lower2.1 ex\hbox{\begin{tikzpicture}[scale=0.3] 
\foreach \x in {1,...,#1}
{ \uuvert{\x} \lvert{\x} }
#2 \end{tikzpicture}}}}
\nc{\disppartnxdd}[2]{{\lower3.6 ex\hbox{\begin{tikzpicture}[scale=0.3] 
\foreach \x in {1,...,#1}
{ \uuuvert{\x} \uuvert{\x} \uvert{\x} \lvert{\x} }
#2 \end{tikzpicture}}}}

\nc{\dispgax}[2]{{\lower0.0 ex\hbox{\begin{tikzpicture}[scale=0.3] 
#2
\foreach \x in {1,...,#1}
{\lvert{\x} }
 \end{tikzpicture}}}}
\nc{\textgax}[2]{{\lower0.4 ex\hbox{\begin{tikzpicture}[scale=0.3] 
#2
\foreach \x in {1,...,#1}
{\lvert{\x} }
 \end{tikzpicture}}}}
\nc{\textlinegraph}[2]{{\raise#1 ex\hbox{\begin{tikzpicture}[scale=0.8] 
#2
 \end{tikzpicture}}}}
\nc{\textlinegraphl}[2]{{\raise#1 ex\hbox{\begin{tikzpicture}[scale=0.8] 
\tikzstyle{vertex}=[circle,draw=black, fill=white, inner sep = 0.07cm]
#2
 \end{tikzpicture}}}}
\nc{\displinegraph}[1]{{\lower0.0 ex\hbox{\begin{tikzpicture}[scale=0.6] 
#1
 \end{tikzpicture}}}}
 
\nc{\disppartnthreeone}[1]{{\lower1.0 ex\hbox{\begin{tikzpicture}[scale=0.3] 
\foreach \x in {1,2,3,5,6}
{ \uvert{\x} }
\foreach \x in {1,2,4,5,6}
{ \lvert{\x} }
\draw[dotted] (3.5,1.5)--(4.5,1.5);
\draw[dotted] (2.5,0)--(3.5,0);
#1 \end{tikzpicture}}}}

\nc{\partn}[4]{\left( \begin{array}{c|c} 
#1 \ & \ #3 \ \ \\ \cline{2-2}
#2 \ & \ #4 \ \
\end{array} \!\!\! \right)}
\nc{\partnlong}[6]{\partn{#1}{#2}{#3,\ #4}{#5,\ #6}} 
\nc{\partnsh}[2]{\left( \begin{array}{c} 
#1 \\
#2 
\end{array} \right)}
\nc{\partncodefz}[3]{\partn{#1}{#2}{#3}{\emptyset}}
\nc{\partndefz}[3]{{\partn{#1}{#2}{\emptyset}{#3}}}
\nc{\partnlast}[2]{\left( \begin{array}{c|c}
#1 \ &  \ #2 \\
#1 \ &  \ #2
\end{array} \right)}

\nc{\uuarcx}[3]{\draw(#1,3)arc(180:270:#3) (#1+#3,3-#3)--(#2-#3,3-#3) (#2-#3,3-#3) arc(270:360:#3);}
\nc{\uuarc}[2]{\uuarcx{#1}{#2}{.4}}
\nc{\uuuarcx}[3]{\draw(#1,4.5)arc(180:270:#3) (#1+#3,4.5-#3)--(#2-#3,4.5-#3) (#2-#3,4.5-#3) arc(270:360:#3);}
\nc{\uuuarc}[2]{\uuuarcx{#1}{#2}{.4}}
\nc{\darcx}[3]{\draw(#1,0)arc(180:90:#3) (#1+#3,#3)--(#2-#3,#3) (#2-#3,#3) arc(90:0:#3);}
\nc{\darc}[2]{\darcx{#1}{#2}{.4}}
\nc{\udarcx}[3]{\draw(#1,1.5)arc(180:90:#3) (#1+#3,1.5+#3)--(#2-#3,1.5+#3) (#2-#3,1.5+#3) arc(90:0:#3);}
\nc{\udarc}[2]{\udarcx{#1}{#2}{.4}}
\nc{\uudarcx}[3]{\draw(#1,3)arc(180:90:#3) (#1+#3,3+#3)--(#2-#3,3+#3) (#2-#3,3+#3) arc(90:0:#3);}
\nc{\uudarc}[2]{\uudarcx{#1}{#2}{.4}}
\nc{\uarcx}[3]{\draw(#1,1.5)arc(180:270:#3) (#1+#3,1.5-#3)--(#2-#3,1.5-#3) (#2-#3,1.5-#3) arc(270:360:#3);}
\nc{\uarc}[2]{\uarcx{#1}{#2}{.4}}
\nc{\darcxhalf}[3]{\draw(#1,0)arc(180:90:#3) (#1+#3,#3)--(#2,#3) ;}
\nc{\darchalf}[2]{\darcxhalf{#1}{#2}{.4}}
\nc{\uarcxhalf}[3]{\draw(#1,1.5)arc(180:270:#3) (#1+#3,1.5-#3)--(#2,1.5-#3) ;}
\nc{\uarchalf}[2]{\uarcxhalf{#1}{#2}{.4}}
\nc{\uarcxhalfr}[3]{\draw (#1+#3,1.5-#3)--(#2-#3,1.5-#3) (#2-#3,1.5-#3) arc(270:360:#3);}
\nc{\uarchalfr}[2]{\uarcxhalfr{#1}{#2}{.4}}

\nc{\bdarcx}[3]{\draw[blue](#1,0)arc(180:90:#3) (#1+#3,#3)--(#2-#3,#3) (#2-#3,#3) arc(90:0:#3);}
\nc{\bdarc}[2]{\darcx{#1}{#2}{.4}}
\nc{\rduarcx}[3]{\draw[red](#1,0)arc(180:270:#3) (#1+#3,0-#3)--(#2-#3,0-#3) (#2-#3,0-#3) arc(270:360:#3);}
\nc{\rduarc}[2]{\uarcx{#1}{#2}{.4}}
\nc{\duarcx}[3]{\draw(#1,0)arc(180:270:#3) (#1+#3,0-#3)--(#2-#3,0-#3) (#2-#3,0-#3) arc(270:360:#3);}
\nc{\duarc}[2]{\uarcx{#1}{#2}{.4}}

\nc{\uv}[1]{\fill (#1,2)circle(.1);}
\nc{\lv}[1]{\fill (#1,0)circle(.1);}
\nc{\stline}[2]{\draw(#1,2)--(#2,0);}
\nc{\tlab}[2]{\draw(#1,2)node[above]{\tiny $#2$};}
\nc{\tudots}[1]{\draw(#1,2)node{$\cdots$};}
\nc{\tldots}[1]{\draw(#1,0)node{$\cdots$};}

\nc{\uvw}[1]{\fill[white] (#1,2)circle(.1);}
\nc{\huv}[1]{\fill (#1,1)circle(.1);}
\nc{\llv}[1]{\fill (#1,-2)circle(.1);}
\nc{\arcup}[2]{
\draw(#1,2)arc(180:270:.4) (#1+.4,1.6)--(#2-.4,1.6) (#2-.4,1.6) arc(270:360:.4);
}
\nc{\harcup}[2]{
\draw(#1,1)arc(180:270:.4) (#1+.4,.6)--(#2-.4,.6) (#2-.4,.6) arc(270:360:.4);
}
\nc{\arcdn}[2]{
\draw(#1,0)arc(180:90:.4) (#1+.4,.4)--(#2-.4,.4) (#2-.4,.4) arc(90:0:.4);
}
\nc{\cve}[2]{
\draw(#1,2) to [out=270,in=90] (#2,0);
}
\nc{\hcve}[2]{
\draw(#1,1) to [out=270,in=90] (#2,0);
}
\nc{\catarc}[3]{
\draw(#1,2)arc(180:270:#3) (#1+#3,2-#3)--(#2-#3,2-#3) (#2-#3,2-#3) arc(270:360:#3);
}

\nc{\arcr}[2]{
\draw[red](#1,0)arc(180:90:.4) (#1+.4,.4)--(#2-.4,.4) (#2-.4,.4) arc(90:0:.4);
}
\nc{\arcb}[2]{
\draw[blue](#1,2-2)arc(180:270:.4) (#1+.4,1.6-2)--(#2-.4,1.6-2) (#2-.4,1.6-2) arc(270:360:.4);
}
\nc{\loopr}[1]{
\draw[blue](#1,-2) edge [out=130,in=50,loop] ();
}
\nc{\loopb}[1]{
\draw[red](#1,-2) edge [out=180+130,in=180+50,loop] ();
}
\nc{\redto}[2]{\draw[red](#1,0)--(#2,0);}
\nc{\bluto}[2]{\draw[blue](#1,0)--(#2,0);}
\nc{\dotto}[2]{\draw[dotted](#1,0)--(#2,0);}
\nc{\lloopr}[2]{\draw[red](#1,0)arc(0:360:#2);}
\nc{\lloopb}[2]{\draw[blue](#1,0)arc(0:360:#2);}
\nc{\rloopr}[2]{\draw[red](#1,0)arc(-180:180:#2);}
\nc{\rloopb}[2]{\draw[blue](#1,0)arc(-180:180:#2);}
\nc{\uloopr}[2]{\draw[red](#1,0)arc(-270:270:#2);}
\nc{\uloopb}[2]{\draw[blue](#1,0)arc(-270:270:#2);}
\nc{\dloopr}[2]{\draw[red](#1,0)arc(-90:270:#2);}
\nc{\dloopb}[2]{\draw[blue](#1,0)arc(-90:270:#2);}
\nc{\llloopr}[2]{\draw[red](#1,0-2)arc(0:360:#2);}
\nc{\llloopb}[2]{\draw[blue](#1,0-2)arc(0:360:#2);}
\nc{\lrloopr}[2]{\draw[red](#1,0-2)arc(-180:180:#2);}
\nc{\lrloopb}[2]{\draw[blue](#1,0-2)arc(-180:180:#2);}
\nc{\ldloopr}[2]{\draw[red](#1,0-2)arc(-270:270:#2);}
\nc{\ldloopb}[2]{\draw[blue](#1,0-2)arc(-270:270:#2);}
\nc{\luloopr}[2]{\draw[red](#1,0-2)arc(-90:270:#2);}
\nc{\luloopb}[2]{\draw[blue](#1,0-2)arc(-90:270:#2);}

\nc{\larcb}[2]{
\draw[blue](#1,0-2)arc(180:90:.4) (#1+.4,.4-2)--(#2-.4,.4-2) (#2-.4,.4-2) arc(90:0:.4);
}
\nc{\larcr}[2]{
\draw[red](#1,2-2-2)arc(180:270:.4) (#1+.4,1.6-2-2)--(#2-.4,1.6-2-2) (#2-.4,1.6-2-2) arc(270:360:.4);
}

\rnc{\H}{\mathrel{\mathscr H}}
\rnc{\L}{\mathrel{\mathscr L}}
\nc{\R}{\mathrel{\mathscr R}}
\nc{\D}{\mathrel{\mathscr D}}
\nc{\J}{\mathrel{\mathscr J}}
\nc{\leqR}{\mathrel{\leq_{\R}}}
\nc{\leqL}{\mathrel{\leq_{\L}}}
\nc{\leqJ}{\mathrel{\leq_{\J}}}

\nc{\ssim}{\mathrel{\raise0.25 ex\hbox{\oalign{$\approx$\crcr\noalign{\kern-0.84 ex}$\approx$}}}}
\nc{\POI}{\mathcal{POI}}
\nc{\wb}{\overline{w}}
\nc{\ub}{\overline{u}}
\nc{\vb}{\overline{v}}
\nc{\fb}{\overline{f}}
\nc{\gb}{\overline{g}}
\nc{\hb}{\overline{h}}
\nc{\pb}{\overline{p}}
\rnc{\sb}{\overline{s}}
\nc{\XR}{\pres{X}{R\,}}
\nc{\YQ}{\pres{Y}{Q}}
\nc{\ZP}{\pres{Z}{P\,}}
\nc{\XRone}{\pres{X_1}{R_1}}
\nc{\XRtwo}{\pres{X_2}{R_2}}
\nc{\XRthree}{\pres{X_1\cup X_2}{R_1\cup R_2\cup R_3}}
\nc{\er}{\eqref}
\nc{\larr}{\mathrel{\hspace{-0.35 ex}>\hspace{-1.1ex}-}\hspace{-0.35 ex}}
\nc{\rarr}{\mathrel{\hspace{-0.35 ex}-\hspace{-0.5ex}-\hspace{-2.3ex}>\hspace{-0.35 ex}}}
\nc{\lrarr}{\mathrel{\hspace{-0.35 ex}>\hspace{-1.1ex}-\hspace{-0.5ex}-\hspace{-2.3ex}>\hspace{-0.35 ex}}}
\nc{\nn}{\tag*{}}
\nc{\epfal}{\tag*{$\Box$}}
\nc{\tagd}[1]{\tag*{(#1)$'$}}
\nc{\ldb}{[\![}
\nc{\rdb}{]\!]}
\nc{\sm}{\setminus}
\nc{\I}{\mathcal I}
\nc{\InSn}{\I_n\setminus\S_n}
\nc{\dom}{\operatorname{dom}} \nc{\codom}{\operatorname{dom}}
\nc{\ojin}{1\leq j<i\leq n}
\nc{\eh}{\widehat{e}}
\nc{\wh}{\widehat{w}}
\nc{\uh}{\widehat{u}}
\nc{\vh}{\widehat{v}}
\nc{\fh}{\widehat{f}}
\nc{\textres}[1]{\text{\emph{#1}}}
\nc{\aand}{\emph{\ and \ }}
\nc{\iif}{\emph{\ if \ }}
\nc{\textlarr}{\ \larr\ }
\nc{\textrarr}{\ \rarr\ }
\nc{\textlrarr}{\ \lrarr\ }

\nc{\comma}{,\ }

\nc{\COMMA}{,\qquad}
\nc{\COMMa}{,\quad }
\nc{\COMma}{,\ \ \ }
\nc{\COmma}{,\ \ }
\nc{\TnSn}{\T_n\setminus\S_n} 
\nc{\TmSm}{\T_m\setminus\S_m} 
\nc{\TXSX}{\T_X\setminus\S_X} 
\rnc{\S}{\mathcal S}

\nc{\T}{\mathcal T} 
\rnc{\P}{\mathcal P} 
\nc{\K}{\mathrel\mathscr K}
\nc{\PB}{\mathcal{PB}} 
\nc{\rank}{\operatorname{rank}}

\nc{\mtt}{\!\!\!\mt\!\!\!}

\nc{\sub}{\subseteq}
\nc{\la}{\langle}
\nc{\ra}{\rangle}
\nc{\mt}{\mapsto}
\nc{\im}{\mathrm{im}}
\nc{\id}{\mathrm{id}}
\nc{\bn}{\mathbf{n}}
\nc{\ba}{\mathbf{a}}
\nc{\bl}{\mathbf{l}}
\nc{\bm}{\mathbf{m}}
\nc{\bk}{\mathbf{k}}
\nc{\br}{\mathbf{r}}
\nc{\al}{\alpha}
\nc{\be}{\beta}
\nc{\ga}{\gamma}
\nc{\Ga}{\Gamma}
\nc{\de}{\delta}
\nc{\ka}{\kappa}
\nc{\lam}{\lambda}
\nc{\Lam}{\Lambda}
\nc{\si}{\sigma}
\nc{\Si}{\Sigma}
\nc{\oijn}{1\leq i<j\leq n}
\nc{\oijm}{1\leq i<j\leq m}

\nc{\comm}{\rightleftharpoons}
\nc{\AND}{\qquad\text{and}\qquad}
\nc{\ANd}{\quad\text{and}\quad}

\nc{\bit}{\begin{itemize}}
\nc{\bitbmc}{\begin{itemize}\begin{multicols}}
\nc{\bmc}{\begin{itemize}\begin{multicols}}
\nc{\emc}{\end{multicols}\end{itemize}}
\nc{\eit}{\end{itemize}}
\nc{\ben}{\begin{enumerate}}
\nc{\een}{\end{enumerate}}
\nc{\eitres}{\end{itemize}}
\nc{\eitmc}{\end{itemize}}

\nc{\set}[2]{\{ {#1} : {#2} \}} 
\nc{\bigset}[2]{\big\{ {#1}: {#2} \big\}} 
\nc{\Bigset}[2]{\left\{ \,{#1} :{#2}\, \right\}}

\nc{\pres}[2]{\la {#1} \,|\, {#2} \ra}
\nc{\bigpres}[2]{\big\la {#1} \,\big|\, {#2} \big\ra}
\nc{\Bigpres}[2]{\Big\la \,{#1}\, \,\Big|\, \,{#2}\, \Big\ra}
\nc{\Biggpres}[2]{\Bigg\la {#1} \,\Bigg|\, {#2} \Bigg\ra}

\newcommand{\pf}{\begin{proof}}
\newcommand{\epf}{\end{proof}}
\nc{\epfres}{\qed}
\nc{\pfnb}{\pf}
\nc{\epfnb}{\bigskip}
\nc{\pfthm}[1]{\bigskip \noindent{\bf Proof of Theorem \ref{#1}}\,\,  } 
\nc{\pfprop}[1]{\bigskip \noindent{\bf Proof of Proposition \ref{#1}}\,\,  } 
\nc{\epfreseq}{\tag*{$\Box$}}

\makeatletter
\newcommand\footnoteref[1]{\protected@xdef\@thefnmark{\ref{#1}}\@footnotemark}
\makeatother

\numberwithin{equation}{section}

\newtheorem{thm}[equation]{Theorem}
\newtheorem{lemma}[equation]{Lemma}
\newtheorem{cor}[equation]{Corollary}
\newtheorem{prop}[equation]{Proposition}

\theoremstyle{definition}

\newtheorem{rem}[equation]{Remark}
\newtheorem{defn}[equation]{Definition}
\newtheorem{eg}[equation]{Example}
\newtheorem{ass}[equation]{Assumption}

\title{Sandwich semigroups in locally small categories I: Foundations}
\author{
Igor Dolinka\footnote{Department of Mathematics and Informatics, University of Novi Sad, Novi Sad, Serbia. {\it Emails:} {\tt dockie\,@\,dmi.uns.ac.rs}, {\tt ivana.djurdjev\,@\,dmi.uns.ac.rs}}, 
Ivana \DJ ur\dj ev,\hspace{-.2em}${}^*$\footnote{Mathematical Institute of the Serbian Academy of Sciences and Arts, Beograd, Serbia. {\it Email:} {\tt djurdjev.ivana\,@\,gmail.com}}\;
James East\footnote{Centre for Research in Mathematics, School of Computing, Engineering and Mathematics, Western Sydney University, Sydney, Australia. {\it Email:} {\tt j.east\,@\,westernsydney.edu.au}}, 
Preeyanuch Honyam\footnote{Department of Mathematics, Chiang Mai University, Chiang Mai, Thailand. {\it Emails:} {\tt preeyanuch.h\,@\,cmu.ac.th}, {\tt kritsada.s\,@\,cmu.ac.th}, {\tt jintana.s\,@\,cmu.ac.th}}, \\ 
Kritsada Sangkhanan,\hspace{-.2em}${}^\mathsection$
Jintana Sanwong,\hspace{-.2em}${}^\mathsection$
Worachead Sommanee\footnote{Department of Mathematics and Statistics, Chiang Mai Rajabhat University, Chiang Mai, Thailand. {\it Email:} {\tt worachead\_som\,@\,cmru.ac.th}}
}

\date{}

\maketitle

\vspace{-0.5cm}

\begin{abstract}
Fix (not necessarily distinct) objects $i$ and $j$ of a locally small category $S$, and write $S_{ij}$ for the set of all morphisms $i\to j$.  Fix a morphism $a\in S_{ji}$, and define an operation $\star_a$ on $S_{ij}$ by $x\star_ay=xay$ for all $x,y\in S_{ij}$.  Then $(S_{ij},\star_a)$ is a semigroup, known as a \emph{sandwich semigroup}, and denoted by~$S_{ij}^a$.  This article develops a general theory of sandwich semigroups in locally small categories.  We begin with structural issues such as regularity, Green's relations and stability, focusing on the relationships between these properties on $S_{ij}^a$ and the whole category $S$.  We then identify a natural condition on~$a$, called \emph{sandwich regularity}, under which the set $\Reg(S_{ij}^a)$ of all regular elements of $S_{ij}^a$ is a subsemigroup of~$S_{ij}^a$.  Under this condition, we carefully analyse the structure of the semigroup $\Reg(S_{ij}^a)$, relating it via pullback products to certain regular subsemigroups of $S_{ii}$ and $S_{jj}$, and to a certain regular sandwich monoid defined on a subset of~$S_{ji}$; among other things, this allows us to also describe the idempotent-generated subsemigroup~$\mathbb E(S_{ij}^a)$ of $S_{ij}^a$.  We also study combinatorial invariants such as the \emph{rank} (minimal size of a generating set) of the semigroups $S_{ij}^a$, $\Reg(S_{ij}^a)$ and $\mathbb E(S_{ij}^a)$; we give lower bounds for these ranks, and in the case of $\Reg(S_{ij}^a)$ and $\mathbb E(S_{ij}^a)$ show that the bounds are sharp under a certain condition we call \emph{MI-domination}.  Applications to concrete categories of transformations and partial transformations are given in Part II.

%


%

{\it Keywords}: Categories, sandwich semigroups, transformation semigroups, rank, idempotent rank.

MSC: 20M50; 18B40; 20M10; 20M17; 20M20; 05E15.
\end{abstract}

\section*{Introduction}\label{sect:intro}

\emph{Sandwich operations} arise in many mathematical contexts: representation theory \cite{Munn1955,GW2015}, classical groups~\cite{Brown1955}, category theory \cite{Muhammed2017_2}, automota theory \cite{Chase1979}, topology \cite{Magill1967_2,MS1978}, computational algebra \cite{EEMP2015}, and more.  They are also foundational in the structure theory of (completely 0-simple) semigroups \cite{Rees1940}.  As a prototypical example, consider two non-empty sets $X,Y$, write $\T_{XY}$ for the set of all functions $X\to Y$, and let $a\in\T_{YX}$ be some fixed function $Y\to X$.  
Then we may define a \emph{sandwich operation}~$\star_a$ on $\T_{XY}$ by $f\star_ag=f\circ a\circ g$, for $f,g\in\T_{XY}$.  This operation is associative, and the resulting \emph{sandwich semigroup} $\TXYa=(\T_{XY},\star_a)$ is one of several examples discussed by Lyapin in his influential monograph \cite[Chapter VII]{Lyapin1960}.  
The semigroups $\TXYa$, and related sandwich semigroups of partial transformations have been studied by several authors over the years; see especially the early work of Magill \cite{Magill1967,MS1975}.  For additional historical background, further motivating examples, and an extensive list of references, the reader is referred to the introduction of \cite{DElinear}.

A further family of examples is furnished by the so-called \emph{semigroup variants}.  A \emph{variant} of a semigroup~$S$ is a semigroup of the form $S^a=(S,\star_a)$, where $a$ is a fixed element of $S$, and  $x\star_ay=xay$, for $x,y\in S$.  The sandwich semigroups $\TXYa$ discussed above are not variants in general, since the underlying set $\T_{XY}$ is only a semigroup if $X=Y$, in which case $\T_{XX}=\T_X$ is the \emph{full transformation semigroup} over $X$: i.e., the semigroup of all functions $X\to X$ under composition.  A general theory of variants has been developed by a number of authors; see especially \cite{Hickey1986,Hickey1983,KL2001}.  For a recent study of variants of finite full transformation semigroups, and for further references and historical discussion, see \cite{DEvariants} and also \cite[Chapter 13]{GMbook}.

The article \cite{DElinear} initiated the study of general sandwich semigroups in arbitrary (locally small) categories.  This allowed many of the situation-specific results previously proved about various families of sandwich semigroups (such as $\TXYa$, mentioned above) to be treated in a unified fashion.  Consider a locally small category $\C$; by this, it is meant that for any objects $X,Y$ of $\C$, the collection $\C_{XY}$ of all morphisms $X\to Y$ forms a set.  (In fact, the results of \cite{DElinear} were formulated for the more general class of \emph{partial semigroups}, in which identities are not required to exist; formal definitions are given in Section \ref{subsect:basic}.)  Let $X,Y$ be two objects of $\C$, and let $a\in\C_{YX}$ be a fixed morphism $Y\to X$.  Then the set $\C_{XY}$ becomes a semigroup under the sandwich operation $\star_a$, defined by $f\star_ag=fag$, for $f,g\in\C_{XY}$.  When $\C$ is the category of non-empty sets and mappings, we obtain the sandwich semigroups $\TXYa$ discussed above.  

The main motivating example in \cite{DElinear} was the category $\M=\M(\F)$ of all finite dimensional matrices over a field $\F$.  Sandwich semigroups in the category $\M$ are closely related to \emph{Munn rings} \cite{Munn1955} and Brown's \emph{generalised matrix algebras} \cite{Brown1955}, both of which are important tools in representation theory.  The article~\cite{DElinear} began by developing some basic theory of arbitrary sandwich semigroups, mostly concerning Green's relations and (von Neumann) regularity; see Section \ref{subsect:Green} below for definitions and more details.  This theory was then applied to the category $\M$ itself, forming the basis for deeper investigations of the linear sandwich semigroups $\MmnA=(\Mmn,\star_A)$; here, $\Mmn=\Mmn(\F)$ denotes the set of all $m\times n$ matrices over $\F$, and~$A\in\Mnm$ is a fixed $n\times m$ matrix.  
Many of these results on the linear sandwich semigroups were combinatorial in nature, taking inspiration from classical results in (linear) transformation semigroup theory \cite{Howie1966,GH1992,Gomes1987,Howie1978,Howie1990,Gray2007,Erdos1967,Dawlings1982}, and included the classification and enumeration of idempotents and regular elements, description of the idempotent-generated subsemigroup, calculation of the minimum number of (idempotent) matrices required to generate certain subsemigroups of $\MmnA$, and so on.  The proofs of these results relied crucially on intimate relationships between linear sandwich semigroups and various classical (non-sandwich) matrix semigroups.  For example, it was shown that the set $\Reg(\MmnA)$ of all regular elements of $\MmnA$ is itself a semigroup that maps onto the monoid $\M_r=\M_r(\F)$ of all $r\times r$ matrices over~$\F$, where $r$ is the rank of $A$; this allowed for a structural description of $\Reg(\MmnA)$ as a kind of ``inflation'' of $\M_r$, with maximal subgroups of $\M_r$ (all isomorphic to general linear groups of various degrees) ``blown up'' into \emph{rectangular groups} in $\Reg(\MmnA)$.


The current pair of articles, this being the first, has two broad goals:
\ben
\item[I.] Here we further develop the general theory of sandwich semigroups in arbitrary (locally small) categories, extending many of the above-mentioned results from \cite{DElinear,DEvariants} to a far more general setting (see below for more details).
\item[II.] In the sequel \cite{Sandwiches2}, we apply the theory developed here to sandwich semigroups in three concrete categories of transformations and partial transformations.  We also show how many well-known families of (non-sandwich) transformation semigroups arise as special cases of the sandwich semigroup construction, so that many new results (and many previously-known ones) concerning these families may be efficiently deduced as corollaries of our general results.  
\een
A future article \cite{DEdiagrams} will apply the general theory to a number of \emph{diagram categories} \cite{Martin2008}, and we hope that the techniques we develop here will prove useful in the investigation of sandwich semigroups in other categories.  We also hope that the very \emph{idea} of a sandwich semigroup will give category theorists new algebraic tools and invariants with which to study their favourite categories.

The current article is organised as follows.
Section \ref{sect:prelim} begins with the basic definitions concerning partial semigroups and sandwich semigroups, as well as some background from \cite{DElinear}, and then gives further preliminary results concerning (von Neumann) regularity, stability and Green's relations on a partial semigroup~$S$, showing how these relate to the corresponding concepts on its sandwich semigroups~$\Sija$; it also contains results showing how certain properties associated to stability and Green's relations are inherited by partial subsemigroups of $S$ under various regularity assumptions.  Section \ref{sect:RegSija} explores the structure of the regular part $\RegSija$ of a sandwich semigroup $\Sija$.  If the sandwich element $a$ satisfies a natural property called \emph{sandwich-regularity} (which is satisfied by all elements in a von Neumann regular category, for example), then $\RegSija$ is a (regular) subsemigroup of $\Sija$, and the structure of $\RegSija$ and the idempotent-generated subsemigroup $\E_a(\Sija)$ of $\Sija$ are governed by intimate relationships with other (sandwich and non-sandwich) semigroups in $S$.  Section \ref{sect:LMC} identifies a natural property, called \emph{MI-domination} (which is satisfied by all our motivating examples in \cite{Sandwiches2}), under which the relationships explored in Section \ref{sect:RegSija} become even more striking; the MI-domination assumption also allows us to give precise formulae for the rank (and idempotent rank, where appropriate) of $\RegSija$ and $\E_a(\Sija)$.  We conclude with two short sections.  Section \ref{sect:inverse} identifies a condition stronger than sandwich-regularity (satisfied by one of our motivating examples in \cite{Sandwiches2}), under which $\RegSija$ is inverse, and in which case the theory developed in Sections \ref{sect:RegSija} and \ref{sect:LMC} simplifies substantially.  Section \ref{sect:rankSija} contains some remarks about the calculation of the rank of a sandwich semigroup $\Sija$.
We work in standard ZFC set theory; see for example \cite[Chapters 1, 5 and 6]{Jech2003}.  The reader may refer to \cite{MacLane1998} for basics on category theory.

\sectiontitle{Preliminaries}\label{sect:prelim}

\subsectiontitle{Basic definitions}\label{subsect:basic}

Recall from \cite{DElinear} that a \emph{partial semigroup} is a $5$-tuple $(S,\cdot,I,\lam,\rho)$ consisting of a class $S$, a partial binary operation $(x,y)\mt x\cdot y$ (defined on some sub-class of $S\times S$), a class $I$, and functions $\lam,\rho:S\to I$, such that, for all $x,y,z\in S$,~
\bit
\item[(i)] $x\cdot y$ is defined if and only if $x\rho=y\lam$,
\item[(ii)] if $x\cdot y$ is defined, then $(x\cdot y)\lam=x\lam$ and $(x\cdot y)\rho=y\rho$,
\item[(iii)] if $x\cdot y$ and $y\cdot z$ are defined, then $(x\cdot y)\cdot z=x\cdot (y\cdot z)$,
\item[(iv)] for all $i,j\in I$, $S_{ij}=\set{x\in S}{x\lam=i,\ x\rho=j}$ is a set.
\eit
For $i\in I$, we write $S_i=S_{ii}$, noting that these sets are closed under $\cdot$, and hence semigroups.  Partial semigroups in which $S$ and $I$ are sets are sometimes called \emph{semigroupoids}; see for example \cite[Appendix B]{Tilson1987}.  We say that a partial semigroup $(S,\cdot,I,\lam,\rho)$ is \emph{monoidal} if in addition to (i)--(iv), 
\bit
\item[(v)] there exists a function $I\to S:i\mt e_i$ such that, for all $x\in S$, $x\cdot e_{x\rho}=x=e_{x\lam}\cdot x$.
\eit
Note that if (v) holds, then $S_i$ is a monoid, with identity $e_i$, for all $i\in I$.  Note also that any semigroup is trivially a partial semigroup (with $|I|=1$).

In the above definition, we allow $S$ and $I$ to be proper classes; although $S$ and $I$ were assumed to be sets in \cite{DElinear}, it was noted that this requirement is not necessary.  If the operation $\cdot$, the class $I$, and the functions $\lam,\rho$ are clear from context, we will often refer to ``the partial semigroup $S$'' instead of ``the partial semigroup $(S,\cdot,I,\lam,\rho)$''.

Recall \cite{DElinear} that any partial semigroup $S$ may be embedded in a monoidal partial semigroup $\Sone$ in an obvious way.  Note that a monoidal partial semigroup is essentially a locally small category (a category in which the morphisms between fixed objects always form a set), described in an Ehresmann-style ``arrows only'' fashion \cite{Ehresmann1965}; thus, we will sometimes use the terms ``monoidal partial semigroup'' and ``(locally small) category'' interchangeably.
If the context is clear, we will usually denote a product $x\cdot y$ by $xy$.

If $x$ is an element of a partial semigroup $S$, then we say $x$ is \emph{(von Neumann) regular} if there exists $y\in S$ such that $x=xyx$; note then that with $z=yxy$, we have $x=xzx$ and $z=zxz$.  So $x$ is regular if and only if the set $V(x)=\set{y\in S}{x=xyx,\ y=yxy}$ is non-empty.  We write $\Reg(S)=\set{x\in S}{x=xyx\ (\exists y\in S)}$ for the class of all regular elements of $S$.  We say $S$ itself is (von Neumann) regular if $S=\Reg(S)$.  Note that there are other meanings of ``regular category'' in the literature \cite{Barr1971,Lawvere2002}, but we use the current definition for continuity with semigroup (and ring) theory; cf.~\cite{FP1971}.


\subsectiontitle{Green's relations}\label{subsect:Green}

Green's relations and preorders play an essential role in the structure theory of semigroups, and this remains true in the setting of partial semigroups.  
Let $S\equiv(S,\cdot,I,\lam,\rho)$ be a partial semigroup.
As in \cite{DElinear}, if $x,y\in S$, then we say
\begin{itemize}\begin{multicols}{2}
\item $x\leq_{\R} y$ if $x=ya$ for some $a\in\Sone$,
\item $x\leq_{\L} y$ if $x=ay$ for some $a\in \Sone$,
\item $x\leq_{\H} y$ if $x\leq_{\R} y$ and $x\leq_{\L} y$,
\item $x\leq_{\J} y$ if $x=ayb$ for some $a,b\in\Sone$.
\end{multicols}\eitmc
If $\gK$ is one of $\R$, $\L$, $\J$ or $\H$, we write ${\gK}={\leq_{\gK}}\cap{\geq_{\gK}}$ for the equivalence relation on $S$ induced by $\leq_{\gK}$.  The equivalence $\D$ is defined by ${\D}={\R}\vee{\L}$; it was noted in \cite{DElinear} that ${\D}={\R}\circ{\L}={\L}\circ{\R}\sub{\J}$.  We also note that ${\leq_{\H}}={\leq_{\R}}\cap{\leq_{\L}}$ and ${\H}={\R}\cap{\L}$.

Recall that for $i,j\in I$, and for a fixed element $a\in S_{ji}$, there is an associative operation $\star_a$ defined on $S_{ij}$ by $x\star_ay=xay$ for all $x,y\in S_{ij}$.  The semigroup $\Sija=(S_{ij},\star_a)$ is called the \emph{sandwich semigroup} of $S_{ij}$ with respect to $a$, and $a$ is called the \emph{sandwich element}.  The article \cite{DElinear} initiated the general study of sandwich semigroups in arbitrary partial semigroups, and the current article takes this study much further.  

We now state a fundamental result from \cite{DElinear} that describes the interplay between Green's relations on a sandwich semigroup $\Sija$ and the corresponding relations on the partial semigroup $S$.
In order to avoid confusion, if $\K$ is any of $\R$, $\L$, $\J$, $\H$ or $\D$, we write $\K^a$ for Green's $\K$-relation on the sandwich semigroup $\Sija$.  For $x\in S_{ij}$, we write 
\[
K_x = \set{y\in S_{ij}}{x\K y} \AND K_x^a = \set{y\in S_{ij}}{x\K^a y}
\]
for the $\K$-class and $\K^a$-class of $x$ in $S_{ij}$, respectively.  
In describing the relationships between $\K$ and $\K^a$, and in many other contexts, a crucial role is played by the sets defined by 
\[
P_1^a = \set{x\in S_{ij}}{xa\R x} \COMMA
P_2^a = \set{x\in S_{ij}}{ax\L x} \COMMA
P_3^a = \set{x\in S_{ij}}{axa\J x} \COMMA
P^a=P_1^a\cap P_2^a.
\]
The next result is \cite[Theorem 2.13]{DElinear}.


\begin{thm}\label{thm:green_Sij}
Let $(S,\cdot,I,\lam,\rho)$ be a partial semigroup, and let $a\in S_{ji}$ where $i,j\in I$.  If $x\in S_{ij}$, then   
\begin{itemize}\begin{multicols}{2}
\itemit{i} $R_x^a = \begin{cases}
R_x\cap P_1^a &\text{if $x\in P_1^a$}\\
\{x\} &\text{if $x\in S_{ij}\sm P_1^a$,}
\end{cases}$
\itemit{ii} $L_x^a = \begin{cases}
L_x\cap P_2^a &\hspace{0.7mm}\text{if $x\in P_2^a$}\\
\{x\} &\hspace{0.7mm}\text{if $x\in S_{ij}\sm P_2^a$,}
\end{cases}
$
\itemit{iii} $H_x^a = \begin{cases}
H_x &\hspace{7.4mm}\text{if $x\in P^a$}\\
\{x\} &\hspace{7.4mm}\text{if $x\in S_{ij}\sm P^a$,}
\end{cases}$
\itemit{iv} $D_x^a = \begin{cases}
D_x\cap P^a &\text{if $x\in P^a$}\\
L_x^a &\text{if $x\in P_2^a\sm P_1^a$}\\
R_x^a &\text{if $x\in P_1^a\sm P_2^a$}\\
\{x\} &\text{if $x\in S_{ij}\sm (P_1^a\cup P_2^a)$,}
\end{cases}$
\itemit{v} $J_x^a = \begin{cases}
J_x\cap P_3^a &\hspace{2.2mm}\text{if $x\in P_3^a$}\\
D_x^a &\hspace{2.2mm}\text{if $x\in S_{ij}\sm P_3^a$.}
\end{cases}$
\end{multicols}\end{itemize}
Further, if $x\in S_{ij}\sm P^a$, then $H_x^a=\{x\}$ is a non-group $\gHa$-class of $\Sija$.  \epfres
\end{thm}

%

Fix some $i,j\in I$.  We say $e\in S_j$ is a \emph{right identity} for $S_{ij}$ if $xe=x$ for all $x\in S_{ij}$.  We say $a\in S_{ji}$ is \emph{right-invertible} if there exists $b\in S_{ij}$ such that $ab$ is a right identity for $S_{ij}$.  \emph{Left identities} and \emph{left-invertible} elements are defined analogously.  


\begin{lemma}\label{lem:al=xi_gen}
Let $S$ be a partial semigroup, and let $a\in S_{ji}$ where $i,j\in I$.  
\bit
\itemit{i} If $a$ is right-invertible, then $P_1^a=S_{ij}$, $P^a=P_2^a$, and ${\R^a}={\R}$ on $\Sija$.
\itemit{ii} If $a$ is left-invertible, then $P_2^a=S_{ij}$, $P^a=P_1^a$, and ${\L^a}={\L}$ on $\Sija$.
\end{itemize}
\end{lemma}

\pf We just prove (i), as a dual/symmetrical argument will give (ii).  Suppose $a$ is right-invertible, and let $b\in S_{ij}$ be as above.  For any $x\in S_{ij}$, we have $x=xab$, giving $x\R xa$, and $x\in P_1^a$.  This gives $P_1^a=S_{ij}$; the second and third claims immediately follow (using Theorem \ref{thm:green_Sij}(i) to establish the third). \epf




\subsectiontitle{Stability and regularity}\label{subsect:stability}

The concept of \emph{stability} played a crucial role in \cite{DElinear}, and will continue to do so here.  Recall that a partial semigroup $S$ is stable if, for all $a,x\in S$,
\[
xa\J x \iff xa\R x \AND ax\J x \iff ax\L x.
\]
In our studies, we require a more refined notion of stability.
We say an element $a$ of a partial semigroup~$S$ is \emph{$\R$-stable} if $xa\J x \implies xa\R x$ for all $x\in S$.  Similarly, we say $a$ is \emph{$\L$-stable} if $ax\J x \implies ax\L x$ for all $x\in S$.  We say $a$ is \emph{stable} if it is both $\R$- and $\L$-stable.

Recall that a semigroup $T$ is \emph{periodic} if for each $x\in T$, some power of $x$ is an idempotent: that is, $x^{2m}=x^m$ for some $m\geq1$.  It is well known that all finite semigroups are periodic; see for example \cite[Proposition 1.2.3]{Howie}.  The proof of the next result is adapted from that of \cite[Theorem A.2.4]{RSbook}, but is included for convenience.  If $a\in S_{ji}$, for some $i,j\in I$, then $S_{ij}a$ is a subsemigroup of $S_i$, and $aS_{ij}$ a subsemigroup of~$S_j$.  We will have more to say about these subsemigroups in Sections \ref{sect:RegSija} and \ref{sect:LMC}.

\begin{lemma}\label{lem:periodic_stable}
Let $S$ be a partial semigroup, and fix $i,j\in I$ and $a\in S_{ji}$.   
\begin{itemize}
\begin{multicols}{2}
\itemit{i} If $aS_{ij}$ is periodic, then $a$ is $\R$-stable.
\itemit{ii} If $S_{ij}a$ is periodic, then $a$ is $\L$-stable.
\end{multicols}\eitmc
\end{lemma}

\pf We just prove (i), as a dual argument will give (ii).  So suppose $aS_{ij}$ is periodic, and that $x\in S$ is such that $xa\J x$.  For $xa$ to be defined, we must have $x\in S_{kj}$ for some $k\in I$.  Evidently, $xa\leqR x$, so we just need to prove the converse.  Since $xa\J x$, one of the following four conditions hold:
\begin{itemize}
\begin{multicols}{2}
\item[(a)] $x=xa$, or
\item[(b)] $x=xav$ for some $v\in S$, or
\item[(c)] $x=uxa$ for some $u\in S$, or
\item[(d)] $x=uxav$ for some $u,v\in S$.
\end{multicols}
\eitmc
Clearly $x\leqR xa$ if (a) or (b) holds.  If (c) holds, then $x=u(uxa)a=u^2(xa)a$, so that (d) holds.  So suppose~(d) holds, and note that we must have $u\in S_k$ and $v\in S_{ij}$.  Since $aS_{ij}$ is periodic, there exists $m\geq1$ such that $(av)^m$ is idempotent.  But then
$x = u^mx(av)^m = u^mx(av)^m(av)^m = x(av)^m \leqR xa$,
as required.~\epf

%
%

%
If $T$ is any semigroup, we write $\Reg(T)=\set{x\in T}{x=xyx\ (\exists y\in T)}$ for the set of all regular elements of~$T$.  Note that $\Reg(T)$ might not be a subsemigroup of $T$.  


\begin{prop}\label{prop:Reg(Sija)}
Let $S$ be a partial semigroup, and fix $i,j\in I$ and $a\in S_{ji}$.  Then $\RegSija \sub P^a \sub P_3^a$.
\begin{itemize}\begin{multicols}{2}
\itemit{i} If $a$ is $\R$-stable, then $P_3^a\sub P_1^a$.
\itemit{ii} If $a$ is $\L$-stable, then $P_3^a\sub P_2^a$.
\itemit{iii} If $a$ is stable, then $P_3^a=P^a$.
\end{multicols}\eitmc
\end{prop}

\pf The assertion $\RegSija \sub P^a \sub P_3^a$ was proved in \cite[Proposition 2.11]{DElinear}.  If $a$ is $\R$-stable, then
\[
x\in P_3^a \implies x\J axa \implies x\leqJ axa \leqJ xa\leqJ x \implies x\J xa \implies x\R xa \implies x\in P_1^a,
\]
giving~(i).  Part~(ii) follows by a dual argument.  If $a$ is stable, then (i) and (ii) give $P_3^a\sub P_1^a\cap P_2^a=P^a$; since we have already noted that $P^a\sub P_3^a$ (for any $a$), this completes the proof of~(iii).~\epf

\begin{cor}\label{cor:JaDa}
Let $S$ be a partial semigroup, and fix $i,j\in I$ and $a\in S_{ji}$.  If $a$ is stable, and if ${\J}={\D}$ in~$S$, then ${\J^a}={\D^a}$ in $\Sija$.
\end{cor}

\pf Let $x\in S_{ij}$.  We must show that $J_x^a=D_x^a$.  This follows immediately from Theorem \ref{thm:green_Sij}(v) if $x\not\in P_3^a$, so suppose $x\in P_3^a$.  Proposition \ref{prop:Reg(Sija)}(iii) gives $P_3^a=P^a$ (since $a$ is stable), and we have $J_x=D_x$ by assumption.  Together with parts (iv) and (v) of Theorem \ref{thm:green_Sij}, it then follows that $J_x^a=J_x\cap P_3^a=D_x\cap P^a=D_x^a$. \epf

We say a semigroup $T$ is \emph{$\R$-stable} (or \emph{$\L$-stable}, or \emph{stable}) if every element of $T$ is $\R$-stable (or $\L$-stable, or stable, respectively).  The proof of \cite[Proposition 2.14]{DElinear} works virtually unmodified to prove the following.

\begin{lemma}\label{lem:aSa_stable}
Let $a\in S_{ji}$.  If each element of $aS_{ij}a$ is $\R$-stable (or $\L$-stable, or stable) in $S$, then $\Sija$ is $\R$-stable (or $\L$-stable, or stable, respectively).  \epfres
\end{lemma}




The next result shows that the sets $P_1^a$, $P_2^a$ and $P^a=P_1^a\cap P_2^a$ are substructures of $\Sija$; it also describes the relationship between $\Reg(S)$ and $\RegSija$.

\begin{prop}\label{prop:Reg(Sija)2}
Let $S$ be a partial semigroup, and fix $i,j\in I$ and $a\in S_{ji}$.  Then 
\begin{itemize}\begin{multicols}{2}
\itemit{i} $P_1^a$ is a left ideal of $\Sija$,
\itemit{ii} $P_2^a$ is a right ideal of $\Sija$,
\itemit{iii} $P^a$ is a subsemigroup of $\Sija$,
\itemit{iv} $\RegSija=P^a\cap\Reg(S)$,
\itemit{v} $\RegSija=P^a \iff P^a\sub\Reg(S)$.
\item[] ~
\end{multicols}\end{itemize}
\end{prop}

\pf (i).  Suppose $x\in P_1^a$ and $u\in S_{ij}$.  Since $x\R xa$, we have $x=xay$ for some $y\in\Sone$.  But then $uax=uaxay$, so that $uax\R uaxa$, giving $u\star_ax=uax \in P_1^a$.  This proves (i), and (ii) follows by duality.

\pfitem{iii}  The intersection of two subsemigroups is always a subsemigroup.

\pfitem{iv}  By Proposition \ref{prop:Reg(Sija)}, we have $\RegSija\sub P^a$, and we clearly have $\RegSija\sub\Reg(S)$.  Conversely, suppose $x\in P^a\cap\Reg(S)$.  Then $x=xyx=uax=xav$, for some $y\in S_{ji}$ and $u,v\in\Sone$.  But then $x=xyx=(xav)y(uax)=x\star_a(vyu)\star_ax$, with $vyu\in S_{ij}$, giving $x\in\RegSija$.

\pfitem{v}  This follows immediately from (iv). \epf

%

\subsectiontitle{Partial subsemigroups}\label{subsect:subs}

Consider a partial semigroup $S\equiv(S,\cdot,I,\lam,\rho)$, and let $T\sub S$.  If $s\cdot t\in T$ for all $s,t\in T$ with $s\rho=t\lam$, then $(T,\cdot,I,\lam,\rho)$ is itself a partial semigroup (where, for simplicity, we have written $\cdot$, $\lam$ and $\rho$ instead of their restrictions to~$T$), and we refer to $T\equiv(T,\cdot,I,\lam,\rho)$ as a \emph{partial subsemigroup} of $S$.  Note that $T_{ij}=T\cap S_{ij}$ for all $i,j\in I$.  Here we give a number of results that show how concepts such as Green's relations and preorders, and $(\R$- and/or~$\L$-) stability may be inherited by partial subsemigroups under certain regularity assumptions.

If $T$ is a partial subsemigroup of $S$, and if $\K$ is one of Green's relations, we will write $\K^S$ and $\K^T$ for Green's $\K$-relations on $S$ and $T$, respectively, and similarly for the preorders $\leq_{\K}$ (in the case that ${\K}\not={\D}$).  The next result adapts and strengthens a classical semigroup fact; see for example \cite[Proposition~2.4.2]{Howie} or \cite[Proposition A.1.16]{RSbook}.


%


\begin{lemma}\label{lem:Tgreen}
Let $T$ be a partial subsemigroup of $S$, let $x,y\in T$, and let $\K$ be any of $\R$, $\L$ or $\H$.
\bit\begin{multicols}{2}
\itemit{i} If $y\in\Reg(T)$, then $x\leq_{\K^S} y \iff x\leq_{\K^T}y$.
\itemit{ii} If $x,y\in\Reg(T)$, then $x\K^S y \iff x\K^Ty$.
\end{multicols}\eit
\end{lemma}

\pf Clearly (ii) follows from (i).  We just prove (i) in the case that ${\K}={\R}$; the case ${\K}={\L}$ is dual, and the case ${\K}={\H}$ follows from the others.  Suppose $x\in T$ and $y\in\Reg(T)$ are such that $x\leq_{\R^S} y$.  
Then $x=ya$ and $y=yzy$ for some $a\in \Sone$ and $z\in T$.  
But then $x=ya=yzya=y(zx)$; since $zx\in T$, it follows that $x\leq_{\R^T} y$.  We have proved that $x\leq_{\R^S} y \implies x\leq_{\R^T}y$; the converse is clear. \epf

\begin{lemma}\label{lem:Tstable}
Let $T$ be a regular partial subsemigroup of $S$, and let $a\in T$.  If $a$ is $\R$-stable (or $\L$-stable, or stable) in $S$, then $a$ is $\R$-stable (or $\L$-stable, or stable, respectively) in $T$.
\end{lemma}

\pf Suppose $a$ is $\R$-stable in $S$.  Let $x\in T$ with $x\rho=a\lam$.  Then
\[
x\J^T xa \implies x\J^S xa \implies x\R^S xa \implies x\R^T xa,
\]
using $\R$-stability of $a$ in $S$, and Lemma \ref{lem:Tgreen}(ii).  The other statements are proved analogously. \epf

In the next result, if $a\in T_{ji}$, where $T$ is a partial subsemigroup of $S$, we write
\[
P_1^a(S)= \set{x\in S_{ij}}{x\R^S xa} \AND P_1^a(T) = \set{x\in T_{ij}}{x\R^T xa}
\]
to distinguish the sets $P_1^a$ on $S$ and $T$, respectively; we use similar notation with respect to $P_2^a$, $P_3^a$ and $P^a$.

\begin{lemma}\label{lem:TP}
Let $T$ be a partial subsemigroup of $S$, and let $a\in T_{ji}$ where $i,j\in I$.  Then
\bit
\itemit{i} $P_1^a(T)\sub P_1^a(S)\cap T$, with equality if $T_{ij}\cup T_{ij}a\sub\Reg(T)$, 
\itemit{ii} $P_2^a(T)\sub P_2^a(S)\cap T$, with equality $T_{ij}\cup aT_{ij}\sub\Reg(T)$, 
\itemit{iii} $P^a(T)\sub P^a(S)\cap T$, with equality $T_{ij}\cup T_{ij}a\cup aT_{ij}\sub\Reg(T)$,
\itemit{iv} $P_3^a(T)\sub P_3^a(S)\cap T$, with equality if $a$ is stable in $S$ and $T_{ij}\cup T_{ij}a\cup aT_{ij}\sub\Reg(T)$.
\eit
\end{lemma}

\pf For (i), $P_1^a(T) = \set{x\in T_{ij}}{x\R^T xa} \sub \set{x\in T_{ij}}{x\R^S xa} = \set{x\in S_{ij}}{x\R^S xa}\cap T = P_1^a(S)\cap T$; if $T_{ij}\cup T_{ij}a\sub\Reg(T)$, then Lemma \ref{lem:Tgreen}(ii) gives equality in the second step.
Part (ii) follows by duality.  Part~(iii) follows from (i) and (ii).  The first assertion in (iv) is again easily checked; the rest follows from~(iii), Proposition~\ref{prop:Reg(Sija)}(iii) and Lemma \ref{lem:Tstable}.~\epf

\begin{rem}
It is easy to check that we also have $P_3^a(T)=P_3^a(S)\cap T$ if ${\J^T}={\J^S}\cap(T\times T)$.  A similar statement may be made regarding $P_1^a(T)$ and $P_2^a(T)$, replacing $\J$ by $\R$ or $\L$, respectively.
\end{rem}



In the next result, if $a\in T_{ji}\sub S_{ji}$, where $T$ is a partial subsemigroup of $S$, and if $\K$ is one of Green's relations, we write ${\K^a}(S)$ and ${\K^a}(T)$ for the $\K^a$-relation on the sandwich semigroups $S_{ij}^a$ and $T_{ij}^a$, respectively.  We use similar notation for ${\K^a}(S)$- and ${\K^a}(T)$-classes, writing $K_x^a(S)$ and $K_x^a(T)$ for $x\in T_{ij}$, as appropriate.

\begin{prop}
Let $T$ be a partial subsemigroup of $S$, and let $a\in T_{ji}$ where $i,j\in I$.  Then
\bit
\itemit{i} ${\R^a}(T)\sub{\R^a}(S)\cap(T\times T)$, with equality if $T_{ij}\cup T_{ij}a\sub\Reg(T)$, 
\itemit{ii} ${\L^a}(T)\sub{\L^a}(S)\cap(T\times T)$, with equality if $T_{ij}\cup aT_{ij}\sub\Reg(T)$, 
\itemit{iii} ${\H^a}(T)\sub{\H^a}(S)\cap(T\times T)$, with equality if $T_{ij}\cup T_{ij}a\cup aT_{ij}\sub\Reg(T)$.
\eit
\end{prop}

\pf We just prove (i), as (ii) will follow by duality, and (iii) from (i) and (ii).  The stated containment is clear.
Now suppose $T_{ij}\cup T_{ij}a\sub\Reg(T)$, and let $x\in T_{ij}$.  We use Theorem \ref{thm:green_Sij} and Lemmas \ref{lem:Tgreen} and~\ref{lem:TP}.  
If $x\not\in P_1^a(T)=P_1^a(S)\cap T$, then
$
R_x^a(T)=\{x\}=\{x\}\cap T=R_x^a(S)\cap T.
$
If $x\in P_1^a(T)$, then $R_x^a(T)  = R_x(T)\cap P_1^a(T) = (R_x(S)\cap T)\cap (P_1^a(S)\cap T) = (R_x(S)\cap P_1^a(S))\cap T = R_x^a(S)\cap T$. \epf

\sectiontitle{Sandwich-regularity and the structure of $\Reg(\Sija)$}\label{sect:RegSija}

This section explores the set $\RegSija$ of regular elements of a sandwich semigroup $\Sija$ in a partial semigroup $S\equiv(S,\cdot,I,\lam,\rho)$.  In particular, we identify a natural condition on the sandwich element $a\in S_{ji}$, called \emph{sandwich-regularity}, under which $\RegSija$ is a subsemigroup of~$\Sija$.  If $a$ is sandwich-regular, and if $b\in V(a)$, then the structure of $\RegSija$ is intimately related to that of a certain regular monoid $W$ contained in the sandwich semigroup $S_{ji}^b$.  This monoid $W$ is simultaneously a homomorphic image and (an isomorphic copy of) a submonoid of $\RegSija$.
There are three main results in this section.  
Theorem \ref{thm:psi} realises $\RegSija$ as a pullback product of certain regular subsemigroups of the (non-sandwich) semigroups $S_i$ and~$S_j$.  
Theorem~\ref{thm:RG} displays $\RegSija$ as a kind of ``inflation'' of the monoid $W$; see Remark~\ref{rem:inflation} for the meaning of this (non-technical) term.
Finally, Theorem \ref{thm:EaEb} describes the idempotent-generated subsemigroup of~$\Sija$ in terms of the idempotent-generated subsemigroup of $W$.  
In the concrete applications we consider in \cite{Sandwiches2},~$W$ is always (isomorphic to) a well-understood ``classical'' monoid.

\subsectiontitle{The basic commutative diagrams}\label{sect:commutative}

Recall that for $a\in S$, we write $V(a)=\set{b\in S}{a=aba,\ b=bab}$ for the (possibly empty) set of inverses of~$a$ in $S$.  Suppose now that $a\in S$ is regular, and fix some $b\in V(a)$.  Let $j=a\lam$ and $i=a\rho$, so that $a\in S_{ji}$ and~$b\in S_{ij}$.  As noted above, $S_{ij}a$ and $aS_{ij}$ are subsemigroups of the (non-sandwich) semigroups $S_i=S_{ii}$ and~$S_j=S_{jj}$, respectively.  It is also the case that $aS_{ij}a$ is a subsemigroup of the sandwich semigroup~$S_{ji}^b$, since, for any $x,y\in S_{ij}$, $(axa)\star_b(aya)=axabaya=a(xay)a$, with $xay\in S_{ij}$.  In fact, $(aS_{ij}a,\star_b)$ is a \emph{monoid} with identity $a$, and its operation $\star_b$ is independent of the choice of $b\in V(a)$; that is, if also $c\in V(a)$, then the operations $\star_b$ and~$\star_c$ are identical on $aS_{ij}a\sub S_{ji}$ (though the full sandwich semigroups $S_{ji}^b=(S_{ji},\star_b)$ and~$S_{ji}^c=(S_{ji},\star_c)$ might be quite different).  To emphasise this independence, we will write $\starb$ for the operation~$\star_b$ on $aS_{ij}a$.  As noted above, for $x,y\in S_{ij}$, we have $axa\starb aya=axaya$ (so that here $a$ becomes the ``bread'' instead of the ``filling''!).  Evidently, the following diagram of semigroup epimorphisms commutes:
\begin{equation}\label{eq:CD}
\includegraphics{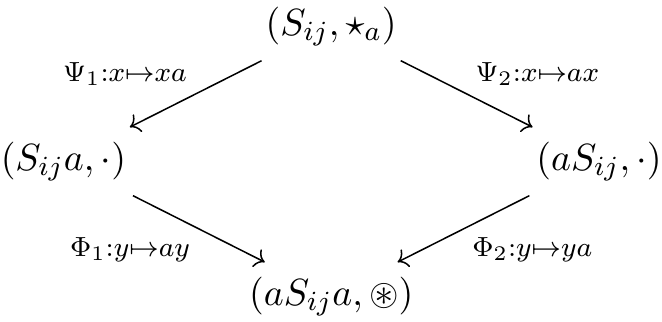}
\end{equation}
Note that $S_{ij}a=(S_{ij}a)ba\sub S_i(ba)=(S_ib)a\sub S_{ij}a$, so that, in fact, $S_{ij}a=S_i(ba)$ is a principal left ideal of~$S_i$: namely, the principal left ideal containing $ba\in S_i$.  Similarly, $aS_{ij}=(ab)S_j$ is the principal right ideal of $S_j$ containing $ab\in S_j$.  Note that the map
\[
\Psi_1:(S_{ij},\star_a)\to(S_{ij}a,\cdot):x\mt xa
\]
is an isomorphism if and only if it is injective, which occurs if and only if $a$ satisfies a certain \emph{right cancellability} law:
\begin{align}
\label{eq:aR}
&(\forall x,y\in S_{ij}) & \hspace{-4cm}xa=ya &\implies x=y.
\intertext{For example, \eqref{eq:aR} holds if $a$ is right-invertible (as defined in Section \ref{subsect:Green}).  Similarly, $\Psi_2$ is an isomorphism if and only if}
\label{eq:aL}
&(\forall x,y\in S_{ij}) & \hspace{-4cm}ax=ay &\implies x=y,
\intertext{which occurs, for example, if $a$ is left-invertible.  As noted by Hickey (in the less general context of \emph{semigroup variants}) in \cite[Lemma 2.2]{Hickey1983}, if $a$ satisfies the weaker cancellativity law}
\label{eq:aRL}
&(\forall x,y\in S_{ij}) & \hspace{-4cm}[xa=ya \text{ and } ax=ay] &\implies x=y,
\end{align}
then the map
\[
\Psi=(\Psi_1,\Psi_2):(S_{ij},\star_a)\to(S_{ij}a,\cdot)\times(aS_{ij},\cdot):x\mt(xa,ax)
\]
is an embedding; although $\Psi$ is not necessarily surjective, its image is a subdirect product (whether or not~$\Psi$ is injective).  Finally, if $a$ satisfies both \eqref{eq:aR} and \eqref{eq:aL}, then the map
\[
\Phi=\Psi_1\Phi_1=\Psi_2\Phi_2:(S_{ij},\star_a)\to(aS_{ij}a,\starb):x\mt axa
\]
is injective, and hence an isomorphism, so that all semigroups in \eqref{eq:CD} are isomorphic: indeed, if $x,y\in S_{ij}$, then $bax,bay\in S_{ij}$, and so 
\[
axa=aya \implies (bax)a=(bay)a \implies bax=bay \implies ax=a(bax)=a(bay)=ay \implies x=y.
\]

In general, even though $a$ is itself regular, none of the semigroups in the diagram \eqref{eq:CD} need be regular semigroups (even if every element of the partial semigroup $S$ is regular).  However, a simple assumption on~$a$ ensures that the set of regular elements in each of these semigroups is a subsemigroup.  We say that the regular element $a\in S_{ji}$ is \emph{sandwich-regular} if $aS_{ij}a\sub\Reg(S)$.  (In particular, if $S$ is regular, then all its elements are sandwich-regular.)


\begin{prop}\label{prop:aPa}
Suppose $i,j\in I$, and that $a\in S_{ji}$ is sandwich-regular.  Let $b\in V(a)$.  Then
\bit
\itemit{i} $\Reg(\Sija)=P^a$ is a (regular) subsemigroup of $\Sija$,
\itemit{ii} $\Reg(S_{ij}a,\cdot)=P^aa=P_2^aa$ is a (regular) subsemigroup of $(S_{ij}a,\cdot)$,
\itemit{iii} $\Reg(aS_{ij},\cdot)=aP^a=aP_1^a$ is a (regular) subsemigroup of $(aS_{ij},\cdot)$, 
\itemit{iv} $(aS_{ij}a,\starb)=aP^aa=aP_1^aa=aP_2^aa$ is a regular subsemigroup of $S_{ji}^b$.
\eit
\end{prop}

\pf (i).  
Let $x\in P^a$.  By sandwich-regularity, we may choose some $y\in V(axa)$.  Since $x\in P^a=P_1^a\cap P_2^a$, we have $x=uax=xav$ for some $u,v\in\Sone$.  But then $x=uaxav=u(axa)y(axa)v=x(aya)x$, so that $x\in\Reg(S)$.  This shows that $P^a\sub\Reg(S)$, and Proposition \ref{prop:Reg(Sija)2}(v) then gives $P^a=\RegSija$.  Proposition~\ref{prop:Reg(Sija)2}(iii) says that $P^a$ is a subsemigroup.

\pfitem{ii}  From (i), it quickly follows that $P^aa\sub\Reg(S_{ij}a,\cdot)$, and that $P^aa$ is a subsemigroup of $(S_{ij}a,\cdot)$.  
If $x\in P_2^a$, then it is easy to show that $xab\in P^a$; since $xa=(xab)a$, this gives $P_2^aa\sub P^aa$.  To complete the proof of (ii), we must show that $\Reg(S_{ij}a,\cdot)\sub P_2^aa$.  So let $g\in\Reg(S_{ij}a,\cdot)$.  Then $g=ghg$ for some $h\in S_{ij}a$, and we may write $g=ya$ and $h=za$ for some $y,z\in S_{ij}$.  Then $g=(yazay)a$, and we have $yazay\in P_2^a$, since $yazay=gzay=(yazaya)zay\leqL a(yazay)\leqL yazay$, which gives $yazay\L a(yazay)$.  

\pfitem{iii}  This is dual to (ii).

\pfitem{iv}  Again, it quickly follows from (i) that $aP^aa$ is a subsemigroup of $S_{ji}^b$; if $x\in P^a$, and if $y\in S_{ij}$ is such that $x=xayax$, then $axa=(axa)\star_b(aya)\star_b(axa)$, so $aP^aa$ is regular.  We also clearly have $aP^aa\sub aP_i^aa\sub aS_{ij}a$, for $i=1,2$.  If $x\in S_{ij}$, then $axa=a(baxab)a$, and it is easy to check that $baxab\in P^a$.  \epf

\begin{rem}\label{rem:baSa}
Although we will not need to know it in this section, it will be convenient to note for later use that the maps
\[
(aS_{ij}a,\starb)\to(baS_{ij}a,\cdot):x\mt bx \AND (aS_{ij}a,\starb)\to(aS_{ij}ab,\cdot):x\mt xb
\]
are isomorphisms, as is easily checked.  We noted above that $S_{ij}a=S_iba$, so it follows that $baS_{ij}a=baS_iba$.  In particular, $(aS_{ij}a,\starb)$ is isomorphic to $(baS_iba,\cdot)$, the \emph{local monoid} of $S_i$ with respect to the idempotent $ba\in S_i$.  (A local monoid of a semigroup $T$, with respect to an idempotent $e$ of $T$, is the monoid $eTe$.)  Similarly, $(aS_{ij}a,\starb)$ is isomorphic to the local monoid $(abS_jab,\cdot)$ of $S_j$.
\end{rem}

{\bf \boldmath For the rest of Section \ref{sect:RegSija}, we fix a sandwich-regular element $a\in S_{ji}$, and an inverse $b\in V(a)$.}  

\medskip

By Proposition \ref{prop:aPa}, the restrictions of the maps in \eqref{eq:CD} yield another commutative diagram of semigroup epimorphisms:
%
\begin{equation}\label{eq:CD_Reg}
\includegraphics{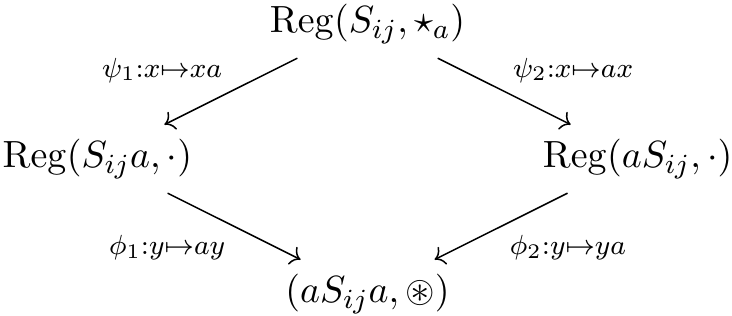}
\end{equation}
The maps in \eqref{eq:CD_Reg} were further explored in \cite{DElinear} in the case that $S=\M(\F)$ was the partial semigroup of (finite dimensional) matrices over a fixed field $\F$; related maps were also studied in \cite{DEvariants}, in the case that $S=\T_X$ was the full transformation semigroup on a finite set $X$.  We wish to explore the diagram \eqref{eq:CD_Reg} here in more generality.  
To simplify notation, for the remainder of Section \ref{sect:RegSija}, we will write
\[
P^a=\Reg(S_{ij},\star_a) \COMMA
W=(aS_{ij}a,\starb)=aP^aa \COMMA
T_1=\Reg(S_{ij}a,\cdot)=P^aa \COMMA
T_2=\Reg(aS_{ij},\cdot)=aP^a.
\]

\subsection[Green's relations on $P^a$ and $W$]{\boldmath Green's relations on $P^a$ and $W$}\label{sect:Green_Pa_W}

In what follows, we will frequently consider Green's relations on the semigroup $P^a=\RegSija$.  We will temporarily denote these by $\K^{P^a}$, where $\K$ is any of $\R,\L,\H,\J,\D$.  Since $P^a$ is regular, we have ${\K^{P^a}}={\K^a}\cap(P^a\times P^a)$ if $\K$ is any of $\R,\L,\H$ (see Lemma \ref{lem:Tgreen}).  The next lemma shows that this is also the case for ${\K}={\D}$; the proof is simple, but included for convenience.


\begin{lemma}\label{lem:PTDJ}
Suppose $T$ is a semigroup for which $Q=\Reg(T)$ is a subsemigroup of $T$.  If $\K$ is any of Green's relations on $T$ other than $\J$, then ${\K^Q}={\K}\cap(Q\times Q)$.
\end{lemma}

\pf This is well known if $\K$ is any of $\R,\L,\H$; see Lemma \ref{lem:Tgreen} (or \cite[Proposition~2.4.2]{Howie} or \cite[Proposition A.1.16]{RSbook}, for example).  Clearly ${\D^Q}\sub{\D}\cap(Q\times Q)$.  Conversely, suppose $(x,y)\in{\D}\cap(Q\times Q)$.  So $x,y\in Q$, and there exists $u\in T$ such that $x\R u\L y$ (in $T$).  Since any element in the $\R$-class (in $T$) of a regular element is regular, it follows that $u\in Q$.  But then $(x,u)\in{\R}\cap(Q\times Q)={\R^Q}$ and $(u,y)\in{\L}\cap(Q\times Q)={\L^Q}$, so that $x\R^Qu\L^Qy$, giving $(x,y)\in{\D^Q}$. \epf

Because of Lemma \ref{lem:PTDJ}, if $\K$ is any of $\R,\L,\H,\D$, then two elements of $P^a=\RegSija$ are $\K^{P^a}$-related in~$P^a$ if and only if they are $\K^a$-related in $\Sija$.  Also, if $x\in P^a$, and if $\K$ is any of $\R,\L,\H,\D$, then $K_x^a\sub P^a$, so that the $\K^{P^a}$-class of $x$ in $P^a$ is precisely the $\K^a$-class of $x$ in $S_{ij}^a$.  As a result, if $\K$ is any of $\R,\L,\H,\D$, we will continue to denote the relations $\K^{P^a}$ on $P^a$ by $\K^a$, and write $K_x^a$ for the $\K^a$-class of $x$ in $P^a$.  In order to avoid confusion, we will write $\JPa$ for Green's $\J$-relation on $P^a$, denote $\JPa$-classes by $J_x^{P^a}$, and so on.


We will also frequently consider Green's relations on the subsemigroup $W=(aS_{ij}a,\starb)=aP^aa$ of $S_{ji}^b$.  These relations can be quite different to those on $S_{ji}^b$ so, in order to avoid confusion, we will denote these by $\K^{\starb}$, where $\K$ denotes any of $\R,\L,\H,\J,\D$, and we denote $\K^\starb$-classes by $K_x^\starb$, and so on.  Although it is not necessarily the case that a $\K^\starb$-class $K_x^\starb$ (in $W$) is equal to the corresponding $\K^b$-class $K_x^b$ (in $S_{ji}^b$), this is at least true for $\H^\starb$-classes, as we now demonstrate.

\begin{lemma}\label{lem:Hb}
If $x\in W=(aS_{ij}a,\starb)$, then $H_x^\starb=H_x^b$.
\end{lemma}

\pf Clearly $H_x^\starb\sub H_x^b$.  Conversely, suppose $y\in H_x^b$.  So $x\R^by$ and $x\L^by$.  So either $x=y$, or else $y=u\star_bx=x\star_bv$ for some $u,v\in S_{ji}$.  Clearly $x=y$ implies $y\in H_x^\starb$, so suppose the latter holds.  Since~$W$ is a monoid with identity $a$, and since $x\in W$, it follows that $y=u\star_bx=u\star_b(x\star_ba)=y\star_ba$, and similarly $y=a\star_by$.  But then $y=a\star_by\star_ba=a(byb)a\in aS_{ij}a=W$.  In particular, $x,y\in W$ are both regular (in~$W$), and so Lemma~\ref{lem:Tgreen}(ii) gives $x\H^\starb y$, whence $y\in H_x^\starb$, as required. \epf


\subsectiontitle{Pullback products and an embedding}\label{sect:pullback}

From the diagram \eqref{eq:CD_Reg}, we may define maps
\[
\phi=\psi_1\phi_1=\psi_2\phi_2:P^a\to W:x\mt axa \AND 
\psi=(\psi_1,\psi_2):P^a\to T_1\times T_2:x\mt(xa,ax).
\]
As a composite of epimorphisms, $\phi$ is itself an epimorphism.  We will study $\phi$ extensively throughout Sections~\ref{sect:RegSija} and \ref{sect:LMC}.
But our first main structural result concerns the map~$\psi$, and shows how it may be used to relate the structure of $P^a$ to that of $T_1$, $T_2$ and $W$.  For the statement, recall that the \emph{pullback product} of semigroups $S_1$ and $S_2$ with respect to a semigroup $T$ and epimorphisms $f_i:S_i\to T$ ($i=1,2$) is the subsemigroup 
$
\set{(s_1,s_2)\in S_1\times S_2}{s_1f_1=s_2f_2}
$
of $S_1\times S_2$.
In particular, pullback products are subdirect.


\begin{thm}\label{thm:psi}
Consider the map $\psi:P^a\to T_1\times T_2:x\mt(xa,ax)$ defined above.  Then
\bit
\itemit{i} $\psi$ is injective,
\itemit{ii} $\im(\psi)=\set{(g,h)\in T_1\times T_2}{ag=ha}=\set{(g,h)\in T_1\times T_2}{g\phi_1=h\phi_2}$.
\eit
In particular, $P^a$ is a pullback product of $T_1$ and $T_2$ with respect to $W$.
\end{thm}

\pf (i).  Suppose $x,y\in P^a$ are such that $(xa,ax)=x\psi=y\psi=(ya,ay)$.  We wish to show that $x=y$.  Since $x,y\in P^a\sub P_1^a$, we have $x=xau$ and $y=yav$ for some $u,v\in S_{ij}$.  Together with $xa=ya$, it follows that $x=yau$ and $y=xav$, whence $x\R^ay$.  A dual argument gives $x\L^ay$.  Green's Lemma in the semigroup~$P^a$ (see for example \cite[Lemma 2.2.1]{Howie}) says that the map $L_x^a\to L_y^a:z\mt zav$ is a bijection, with inverse $L_y^a\to L_x^a:z\mt zau$; so $z=zavau$ for all $z\in L_x^a$.  But $x,y\in L_x^a$ (as $x\L^ay$); thus, combined with $xa=ya$, we obtain $x=xavau=yavau=y$.

%

\pfitem{ii}  First note that any element $(g,h)=(xa,ax)\in\im(\psi)$ clearly satisfies $ag=ha$.  Conversely, suppose $g\in T_1$ and $h\in T_2$ are such that $ag=ha$.  Since $W$ is regular and $ag\in W\sub S_{ji}$, we may choose some $w\in S_{ij}$ with $ag=agwag$.  Since $g\in T_1=\im(\psi_1)$, we may choose some $u,v\in P^a=\Reg(\Sija)$ such that $g=ua$ and $u=uavau$.  Then
\[
g = ua = uavaua = uavag = uavagwag = gwag.
\]
A similar calculation, bearing in mind that $ha=hawha$ (as $ha=ag$), shows that $h=hawh$.  Now put $x=gwh$.  Then $(g,h)=(gwag,hawh)=(gwha,agwh)=(xa,ax)$.  But $x\in P^a=\Reg(\Sija)$, since
\[
x\star_a w\star_a x = (gwh)awa(gwh) = (gwag)w(hawh)=gwh=x,
\]
so $(g,h)=x\psi$ with $x\in P^a$, as required. \epf

\subsectiontitle{Rectangular groups and the internal structure of the $\D^a$-classes of $P^a$}\label{sect:Gh}

We now wish to investigate the way in which the structures of $P^a=\Reg(\Sija)$ and $W=\Reg(aS_{ij}a,\starb)$ may be related via the epimorphism $\phi:P^a\to W:x\mt axa$.  Our next main structural result is Theorem \ref{thm:RG}, which shows that $P^a$ is a kind of ``inflation'' of $W$, in a non-technical sense.  Here we mean that the $\D$- and $\J$-relations on $P^a$ are equivalent to the corresponding relations on $W$; that $\R^\starb$-, $\L^\starb$- and $\H^\starb$-classes of~$W$ yield multiple $\R^a$-, $\L^a$- and $\H^a$-classes of $P^a$; that group/non-group $\H^\starb$-classes of $W$ correspond to group/non-group $\H^a$-classes of $P^a$; and that group $\H^\starb$-classes of $W$ are inflated into rectangular groups in~$P^a$.  For more details, see Theorem \ref{thm:RG}, Remark \ref{rem:inflation} and Figure \ref{fig:inflation}.

For simplicity, if $x\in P^a$, we will write
\[
\xb=x\phi=axa\in W.
\]
We also extend this notation to subsets of $P^a$: if $X\sub P^a$, then we write $\Xb=\set{\xb}{x\in X}$.

If $\K$ is any of $\R,\L,\J,\H,\D$, and if $x,y\in P^a$, we will write $x\gKh^a y$ if $\xb \K^\starb \yb$ in $W$.  For $x\in P^a$, we will write $\Kh_x^a$ for the $\gKh^a$-class of $x$ in~$P^a$.  Since $\Kh_x^a=K_x^\starb \phi^{-1}$ is the preimage of a $\K^\starb $-class of $W$, it follows that $\Kh_x^a$ is a union of $\K^a$-classes of~$P^a$.  The next result demonstrates a number of relationships between the relations $\K^a$ and $\gKh^a$.  

For $X\sub P^a$ and $Y\sub W$, we write $E_a(X)=\set{x\in X}{x=xax}$ and $E_b(Y)=\set{y\in Y}{y=yby}$ for the set of $\star_a$-idempotents from $X$ and the set of $\star_b$-idempotents from $Y$, respectively.
For $u\in P^a$, we write 
\[
V_a(u)=\set{v\in P^a}{u=u\star_av\star_au,\ v=v\star_au\star_av}\]
for the set of all $\star_a$-inverses of $u$ from~$P^a$. 




\begin{lemma}\label{lem:K_hat}
We have
\begin{itemize}\begin{multicols}{2}
\itemit{i} ${\R^a}\sub{\gRh^a}\sub{\D^a}$,
\itemit{ii} ${\L^a}\sub{\gLh^a}\sub{\D^a}$,
\itemit{iii} ${\H^a}\sub{\gHh^a}\sub{\D^a}$, 
\itemit{iv} ${\gDh^a}={\D^a}\sub{\gJh^a}={\JPa}$.
\end{multicols}\end{itemize}
\end{lemma}

\pf (i).  First, it is clear that ${\R^a}\sub{\gRh^a}$.  Now let $(x,y)\in{\gRh^a}$.  Since $P^a$ is regular, $x\R^a e$ and $y\R^a f$ for some $e,f\in E_a(P^a)$.  Since ${\R^a}\sub{\D^a}$, it is enough to show that $e\D^af$.  To do this, we will show that $e\R^a eaf \L^a f$.  Evidently, $eaf\leqRa e$ and $eaf\leqLa f$, so it remains to show the reverse inequalities.  Now, $e\gRh^ax\gRh^ay\gRh^af$, so $\eb\R^\starb \fb$ (in $W$).  Since $\eb,\fb\in E_b(W)$, it follows that $\eb=\fb\starb\eb$ and $\fb=\eb\starb\fb$: that is, $aea=afaea$ and $afa=aeafa$.  But then
\[
e=eaeae=eafaeae\leqRa eaf \AND f=fafaf=faeafaf=faeaf\leqLa eaf,
\]
as required.  

%

\pfitem{ii) and (iii}   First, (ii) follows from (i) by duality, and then (iii) follows from (i) and (ii), since ${\gHh^a}={\gRh^a}\cap{\gLh^a}$.

\pfitem{iv}  It is clear that ${\D^a}\sub{\gDh^a}\sub{\gJh^a}$ and that ${\JPa}\sub{\gJh^a}$, so it remains to show that ${\gDh^a}\sub{\D^a}$ and ${\gJh^a}\sub{\JPa}$.  The first of these follows from (i) and (ii), since ${\gDh^a}={\gRh^a}\vee{\gLh^a}$.  For the second, suppose $(x,y)\in{\gJh^a}$.  
So $\xb\J^\starb \yb$.  In particular, $\xb\leq_{\J^\starb }\yb$, so there exist $u,v\in P^a$ with $\xb=\ub\starb\yb\starb\vb$: that is, $axa=auayava$.  
%
Choose some $z\in V_a(x)$.  Then
\[
x=xazax=xazaxazax=xaz(auayava)zax\leqJPa y.
\]
Similarly, $\yb\leq_{\J^\starb }\xb$ implies $y\leqJPa x$, so that $(x,y)\in{\JPa}$. \epf

As a consequence of the above proof of Lemma \ref{lem:K_hat}(iv), we also have the following.

\begin{cor}
If $x,y\in P^a$, then $J_x^{P^a}\leqJPa J_y^{P^a} \iff J_{\xb}^\starb \leq_{\J^\starb } J_{\yb}^\starb$. 
\epfres
\end{cor}

The next result will be of fundamental importance in what follows.  Since idempotents are regular, we clearly have $E_a(P^a)=E_a(\Sija)$.

\begin{lemma}\label{lem:EaEb}
We have $E_a(P^a)=E_a(\Sija)=E_b(W)\phi^{-1}$.
\end{lemma}

\pf Clearly $E_a(P^a)\sub E_b(W)\phi^{-1}$.  Conversely, suppose $x\in E_b(W)\phi^{-1}$, so $x\in P^a$ and $axa=(axa)b(axa)=axaxa$.  Then, for any $y\in V_a(x)$,
$x=xayax = xa(yaxay)ax=xay(axaxa)yax=xax$.~\epf


We are now ready to prove our second main structural result concerning $P^a=\RegSija$.  Recall that a semigroup $T$ is a \emph{rectangular band} if it is isomorphic to a semigroup with underlying set $I\times J$ (where $I,J$ are non-empty sets), and with product defined by $(i_1,j_1)(i_2,j_2)=(i_1,j_2)$; if $|I|=r$ and $|J|=l$, we say that~$T$ is an $r\times l$ rectangular band.  It is well known that $T$ is a rectangular band if and only if $x=x^2=xyx$ for all $x,y\in T$; see for example \cite[Page 7]{Howie}.
A \emph{rectangular group} is (isomorphic to) the direct product of a rectangular band and a group.  By the \emph{dimensions} of a rectangular band or group, we mean $r=|I|$ and~$l=|J|$, as above; these are the number of $\R$- and $\L$-classes, respectively.
Recall that $H_x^\starb=H_x^b$ for all~$x\in W$, by Lemma~\ref{lem:Hb}.

\begin{thm}\label{thm:RG}
Let $x\in P^a$, and put $r=|\Hh_x^a/{\R^a}|$ and $l=|\Hh_x^a/{\L^a}|$.  Then
\bit
\itemit{i} the restriction to $H_x^a$ of the map $\phi:P^a\to W$ is a bijection $\phi|_{H_x^a}:H_x^a\to H_{\xb}^b$,
\itemit{ii} $H_x^a$ is a group if and only if $H_{\xb}^b$ is a group, in which case these groups are isomorphic, 
\itemit{iii} if $H_x^a$ is a group, then $\Hh_x^a$ is an $r\times l$ rectangular group over $H_x^a$, 
\itemit{iv} if $H_x^a$ is a group, then $E_a(\Hh_x^a)$ is an $r\times l$ rectangular band.
\eit
\end{thm}

\pf We begin with (ii).  First, it is clear that $H_x^a$ being a group implies $H_{\xb}^b$ is a group.  Conversely, suppose~$H_{\xb}^b$ is a group, and let its identity be $\yb$, where $y\in P^a$.  Let the inverse of $\xb$ in $H_{\xb}^b$ be $\zb$, where~$z\in P^a$.  Put
\[
w = xazayazax.
\]
We claim that $w\in H_x^a$ and that $w\in E_a(P^a)$, from which it will follow that $H_x^a$ is a group.  Indeed, first note that
\[
\wb = \xb \starb \zb \starb \yb \starb \zb \starb \xb = \yb \starb \yb \starb \yb = \yb,
\]
as $\zb$ and $\xb$ are inverses in $H_{\xb}^b$, and as $\yb$ is the identity of this group.  It follows from Lemma \ref{lem:EaEb} that $w\in E_a(P^a)$.  To show that $w\in H_x^a$, we must show that $w\R^a x$ and $w\L^a x$.  We show just the first of these, and the second will follow by a symmetrical argument.  Evidently, $w\leqRa x$.  To prove the converse, choose some $u\in V_a(x)$, and note that
\begin{align*}
x = xauax &= xauaxauax = xau\cdot \xb \cdot uax = xau\cdot \yb\starb\xb \cdot uax = xau\cdot \wb\starb\xb \cdot uax \\
&= xau\cdot (axazayazaxa)b(axa) \cdot uax = xazayazaxax = wax,
\end{align*}
so that $x\leqRa w$, as required.  As noted above, this completes the proof that $H_x^a$ is a group (with identity~$w$).
It now follows that the restriction to $H_x^a$ of the map $\phi:P^a\to W$ is a group homomorphism $\phi|_{H_x^a}:H_x^a\to H_{\xb}^b$.  By Lemma~\ref{lem:EaEb}, the group theoretic kernel of $\phi|_{H_x^a}$---that is, the set $\wb\phi^{-1}=\set{u\in H_x^a}{\ub=\wb}$---is equal to~$\{w\}$.  It follows from standard group theoretical facts that $\phi|_{H_x^a}$ is injective.  To complete the proof of (ii), it remains to show that $\phi|_{H_x^a}$ is surjective.  With this in mind, let $q\in H_{\xb}^b$, and write $q=\ub$, where $u\in P^a$.  Put 
\[
v=wauaw.
\]
Since $\vb=\wb\starb\ub\starb\wb=\ub=q$, it remains only to show that $v\in H_x^a$.  Since $H_x^a=H_w^a$, it suffices to show that $v\R^a w$ and $v\L^a w$.  Again, we just do the former.  Evidently $v\leqRa w$.  Let $\sb$ be the inverse of $\ub$ in $H_{\xb}^b$, where $s\in P^a$.  Then
\[
w = wawaw = w\cdot \wb \cdot w = w\cdot \ub\starb\sb \cdot w = w\cdot \ub\starb\wb\starb\sb \cdot w = wauawasaw = vasaw \leqRa v,
\]
as required.  This completes the proof of (ii).

\pfitem{i}  Note that we no longer assume $H_x^a$ is a group.  Since $P^a$ is regular, we may choose some idempotent $y\in E_a(P^a)$ such that $x\R^a y$, noting that $\phi|_{H_y^a}:H_y^a\to H_{\yb}^b$ is a bijection, by part (ii).  Choose some $u,v\in P^a$ such that $y=x\star_a u$ and $x=y\star_a v$.  By Green's Lemma in the semigroups $P^a$ and $W$ (see \cite[Lemma~2.2.1]{Howie}), the maps
\begin{align*}
&\th_1:H_x^a\to H_y^a:z\mt z\star_a u ,&
&\th_3:H_{\xb}^b\to H_{\yb}^b:\zb\mt \zb\starb \ub ,\\
&\th_2:H_y^a\to H_x^a:z\mt z\star_a v, &
&\th_4:H_{\yb}^b\to H_{\xb}^b:\zb\mt \zb\starb \vb 
\end{align*}
are bijections, with $\th_2=\th_1^{-1}$ and $\th_4=\th_3^{-1}$.  Then $\phi|_{H_x^a} = \th_1 \circ\phi|_{H_y^a} \circ\th_4$, a composite of bijections, since for any $q\in H_x^a$,
\[
q \lmap{\th_1} qau \lmap{\phi} aqaua \lmap{\th_4} aqauabava = aqauava = a(q\th_1\th_2)a = aqa = q\phi.
\]  
%
%

\pfitem{iii) and (iv}  Suppose $H_x^a$ is a group.  Without loss of generality, we may assume that $x$ is the identity of this group, so also $\xb$ is the identity of $H_{\xb}^b$.  For convenience, we will write $T=\Hh_x^a$.  To show that $T$ is a rectangular group, it suffices to show that:
\bit
\item[(a)] $T$ is a semigroup,
\item[(b)] $T$ is a union of groups,
\item[(c)] the set $E_a(T)$ of idempotents of $T$ forms a rectangular band (which of course proves (iv)).
\eit
Now, $T=H_{\xb}^b\phi^{-1}$ is a union of $\H^a$-classes, and we know that each of these are groups, by (ii), so (b) certainly holds.  To show (a), suppose $u,v\in T$.  Then $\overline{u\star_a v}=\ub\starb\vb\in H_{\xb}^b$, since $H_{\xb}^b$ is a group, so it follows that $u\star_av\gHh^a x$, whence $u\star_av\in T$.  For (c), let $y,z\in E_a(T)$.  As noted before the statement of the result, it is enough to show that $y\star_az\star_ay=y$.  Now, $\yb=\zb=\xb$, as $\xb$ is the unique idempotent of $H_{\xb}^b$, so 
$\yb\starb\zb\starb\yb=\xb=\yb$.  But then
\[
y=yayay=y\cdot\yb\cdot y = y\cdot\yb\starb\zb\starb\yb\cdot y= y\cdot ayazaya\cdot y=yazay=y\star_az\star_ay,
\]
as required.  \epf

\begin{rem}\label{rem:inflation}
By the preceding series of results, the structure of $P^a=\RegSija$, in terms of Green's relations, is a kind of ``inflation'' of the corresponding structure of $W=\Reg(aS_{ij}a,\starb)$:
\bit
\item The partially ordered sets $(P^a/{\JPa},\leqJPa)$ and $(W/{\J^\starb },\leq_{\J^\starb })$ are order-isomorphic.
\item The sets $P^a/{\D^a}$ and $W/{\D^\starb }$ are in one-one correspondence, via $D_x^a\mt D_{\xb}^\starb $.
\item Each $\gKh^a$-class in $P^a$ is a union of $\gK^a$ classes.
\item The $\R^\starb $-, $\L^\starb $- and $\H^\starb $-classes contained within a single $\D^\starb $-class $D_{\xb}^\starb $ of $W$ are in one-one correspondence with the $\gRh^a$-, $\gLh^a$- and $\gHh^a$-classes in the ${\gDh^a}=\D^a$-class $D_x^a$ of $P^a$. 
\item An $\gHh^a$-class $\Hh_x^a$ is a union of $\H^a$-classes, and these are either all non-groups (if $H_{\xb}^\starb=H_{\xb}^b$ is a non-group $\H^\starb $-class of $W$) or else all groups (if $H_{\xb}^\starb=H_{\xb}^b$ is a group); in the latter case, $\Hh_x^a$ is a rectangular group.
\eit
Figure \ref{fig:inflation} illustrates the last two points in a so-called \emph{egg-box diagram}.  In the left egg-box, which displays a ${\gDh^a}={\D^a}$-class in $P^a$, $\R^a$-related elements are in the same row, $\L^a$-related elements in the same column, and $\H^a$-related elements in the same cell; a similar convention is used in the right egg-box, which displays the corresponding $\D^\starb$-class in $W$.  Group $\H^a$- and $\H^\starb $-classes are shaded gray, and solid lines in the left egg-box denote boundaries between ${\gRh^a}$-classes and ${\gLh^a}$-classes.  See also \cite[Figures 2--7]{Sandwiches2}, which display egg-box diagrams of entire sandwich semigroups in the \emph{partial transformation category} $\PT$.
\end{rem}

\begin{figure}[ht]
\begin{center}
\scalebox{.8}{
\begin{tikzpicture}[scale=1]
\node (D1) at (0,0) {\DaClass{12}{8}{
1/3,1/4,1/5,1/6,1/7,1/9,1/10,1/11,1/12,
2/3,2/4,2/5,2/6,2/7,2/9,2/10,2/11,2/12,
3/3,3/4,3/5,3/6,3/7,3/9,3/10,3/11,3/12,
4/3,4/4,4/5,4/6,4/7,4/9,4/10,4/11,4/12,
5/8,
6/1,6/2,6/3,6/4,6/9,6/10,6/11,6/12,
7/1,7/2,7/3,7/4,7/9,7/10,7/11,7/12,
8/1,8/2,8/3,8/4,8/9,8/10,8/11,8/12
}{}{1}
{0,2,4,7,8,12}
{0,3,4,8}
};
\node (D2) at (13,0) {\DClass{5}{3}{1/2,1/3,1/5,2/4,3/1,3/2,3/5}{}{1}};
\end{tikzpicture}
}
\end{center}
\vspace{-5mm}
\caption{A ${\D^a}$-class of $P^a$ (left) and its corresponding $\D^\starb $-class of $W$ (right).  See Remark \ref{rem:inflation} for more information.}
\label{fig:inflation}
\end{figure}
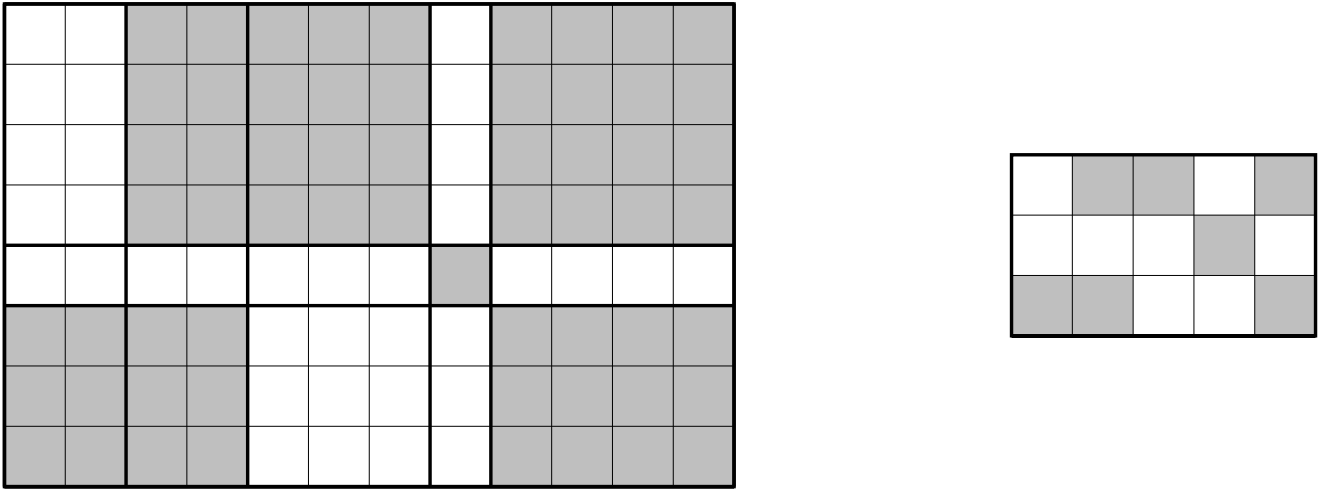

\subsectiontitle{Generation and idempotent-generation}\label{sect:gen}

Recall that Lemma \ref{lem:EaEb} says that the preimage under $\phi:P^a\to W$ of any idempotent of $W$ is an idempotent of $P^a$.  Our next main result, the last of Section \ref{sect:RegSija}, extends this to \emph{products} of idempotents; see Theorem~\ref{thm:EaEb}.  
The key step in the proof of Theorem \ref{thm:EaEb} is the next technical lemma, which will also be useful in other contexts.
For $X\sub S_{ij}$, we write $\la X\ra_a$ for the $\star_a$-subsemigroup of $\Sija$ generated by $X$ (and use similar notation for subsemigroups of $S_{ji}^b$).

\begin{lemma}\label{lem:Xb}
If $X\sub P^a$, then $\la\Xb\ra_b\phi^{-1}\sub\la X\cup E_a(P^a)\ra_a$.
\end{lemma}

\pf Suppose $x\in\la\Xb\ra_b\phi^{-1}$.  So $x\in P^a$ and $\xb\in\la\Xb\ra_b$.  We may therefore write $\xb=\xb_1\starb\cdots\starb\xb_k$, where $x_1,\ldots,x_k\in X$.  Put $y=x_1\star_a\cdots\star_ax_k$.  To prove that $x\in\la X\cup E_a(P^a)\ra_a$, it suffices to show that $x=p\star_ay\star_av$ for some $p,v\in E_a(P^a)$.  

Since $\yb=\xb$, we certainly have $\xb\H^\starb \yb$, so that $x\gHh^a y$.  Since ${\gHh^a}\sub{\D^a}$, by Lemma~\ref{lem:K_hat}(iii), we also have $x\D^ay$.  Since $P^a$ is regular, we may choose some idempotents $u,v\in E_a(P^a)$ such that $u\R^ay$ and $v\L^ax$.  
Since $x\D^ay\R^au$, it follows that $x\D^au$.
Consequently, $R_x^a\cap L_u^a$ is an $\H^a$-class in $D_x^a$: say, $R_x^a\cap L_u^a=H_p^a$.  Since $p\L^au$, and since ${\L^a}\sub{\gLh^a}$, it follows that $p\gLh^au$.  Also, $p\R^ax\gRh^ay\R^au$, so $p\gRh^au$.  Thus, $p\in\Hh_u^a$.  Since $u\in E_a(P^a)$, Theorem~\ref{thm:RG}(iii) tells us that $\Hh_u^a$ is a rectangular group.  In particular, $H_p^a\sub\Hh_u^a$ is itself a group.  Thus, without loss of generality, we may assume that $p$ is the identity of this group; in particular, $p\in E_a(P^a)$.  
By a dual argument, we may show that $L_y^a\cap R_v^a=H_q^a$, for some idempotent $q\in E_a(P^a)$.  

Since $p\L^au$ and $q\R^av$, with $u,q\in E_a(P^a)$, we have $p=p\star_au$ and $v=q\star_av$.  Thus, by Green's Lemma in the semigroup~$P^a$ (see \cite[Lemma 2.2.1]{Howie}), the maps
\[
\th_1:R_u^a\to R_p^a:z\mt p\star_az \AND \th_2:L_q^a\to L_v^a:z\mt z\star_av
\]
are bijections; in addition, $(z,z\th_1)\in{\L^a}$ for all $z\in R_u^a$, while $(z,z\th_2)\in{\R^a}$ for all $z\in L_q^a$.
Now, $y\in R_u^a$, so $p\star_ay=y\th_1\in R_p^a\cap L_y^a=R_x^a\cap L_q^a$.  Since $p\star_ay\in L_q^a$, we also have 
\[
p\star_ay\star_av=(p\star_ay)\th_2\in L_v^a\cap R_{p\star_ay}^a=L_x^a\cap R_x^a=H_x^a.
\]
Since $p\gHh^au$, and since $p,u\in E_a(P^a)$, it follows that $\pb=\ub$ (both being equal to the identity element of~$H_{\ub}^b$).  Similarly, $\qb=\vb$.  Since $y\R^au$ and $y\L^aq$, with $u,q\in E_a(P^a)$, we also have $y=u\star_ay=y\star_aq$, from which it follows that $\yb=\ub\starb\yb=\yb\starb\qb$.  Putting this all together, we obtain
\[
\overline{p\star_ay\star_av} = \pb\starb\yb\starb\vb = \ub\starb\yb\starb\qb = \yb = \xb.
\]
That is, $(p\star_ay\star_av)\phi=x\phi$.  But $\phi$ is injective on $H_x^a$, by Theorem~\ref{thm:RG}(i), and we showed earlier that $p\star_ay\star_av\in H_x^a$, so it follows that $x=p\star_ay\star_av$, as required. \epf

As noted above, one important application of Lemma \ref{lem:Xb} is to describe the idempotent-generated subsemigroup of~$P^a$ in terms of the idempotent-generated subsemigroup of $W$; the last result of this section gives this description.  We write $\E_a(P^a)=\la E_a(P^a)\ra_a$ and ${\E_b(W)=\la E_b(W)\ra_b}$ for these idempotent-generated subsemigroups.  Again, it is clear that~${\E_a(\Sija)=\E_a(P^a)}$.

\begin{thm}\label{thm:EaEb}
We have $\E_a(\Sija)=\E_a(P^a)=\E_b(W)\phi^{-1}$.
\end{thm}

\pf Clearly $\E_a(P^a)\sub\E_b(W)\phi^{-1}$.  The converse follows quickly from Lemma \ref{lem:Xb} with $X=E_a(P^a)$, keeping in mind that $\overline{E_a(P^a)}=E_b(W)$, by Lemma \ref{lem:EaEb}. \epf

\sectiontitle{MI-domination and the ranks of $\Reg(\Sija)$ and $\E_a(\Sija)$}\label{sect:LMC}

In Section \ref{sect:RegSija}, we gave a thorough structural description of the regular and idempotent-generated subsemigroups, $\RegSija$ and $\E_a(\Sija)$, of the sandwich semigroup~$\Sija$, under the assumption that the sandwich element $a\in S_{ji}$ was sandwich-regular.  
%
%
The main purpose of the current section is to prove results concerning the rank (and idempotent rank, where applicable) of these subsemigroups.  The main results---see Theorems~\ref{thm:rankU} and \ref{thm:rankE}---give lower bounds for these ranks, in terms of the (idempotent) ranks of~${W=(aS_{ij}a,\starb)}$ and its idempotent-generated subsemigroup $\E_b(W)$, as well as other parameters such as the dimensions of the rectangular group $\Hh_b^a$.  
The lower bounds turn out to be exact values in the case that $\RegSija$ satisfies a natural condition we formulate below and call \emph{MI-domination}; this is a natural extension of the so-called \emph{RP-domination} property from \cite{BH1984}.  In all of the motivating examples we study in \cite{Sandwiches2}, $\RegSija$ is always MI-dominated; for natural examples where MI-domination fails to hold, see \cite{DEdiagrams}.

\medskip

{\bf \boldmath For the rest of Section \ref{sect:LMC}, we fix a sandwich-regular element $a\in S_{ji}$, and an inverse $b\in V(a)$.}

\subsectiontitle{MI-domination}

We begin by defining the MI-domination property, and discussing some related concepts.
Let $T$ be a regular semigroup, and let $u\in T$.  Recall that $u$ is \emph{regularity-preserving} if the \emph{variant semigroup} $T^u=(T,\star_u)$ is regular \cite{Hickey1983}, and that $u$ is a \emph{mid-identity} if $xuy=xy$ for all $x,y\in T$ \cite{Yamada1955}.  We write
\[
\RP(T) = \set{u\in T}{\text{$u$ is regularity-preserving}} \AND \MI(T) = \set{u\in T}{\text{$u$ is a mid-identity}}.
\]
Clearly $\MI(T)\sub\RP(T)$, since if $u\in\MI(T)$, then the operation $\star_u$ is precisely the original operation of~$T$.  It is easy to see that $\MI(T)$ is a rectangular band; that is, $u=u^2=uvu$ for all $u,v\in\MI(T)$.  Indeed, $u^2=uvu$ is clear, and if $u\in\MI(T)$, then for any $x\in V(u)$, we have $u=(ux)u=(uux)u=uu$.
In particular, $\MI(T)\sub E(T)$.
As explained in \cite{KL2001,Hickey1983}, the motivation for studying regularity-preserving elements is partly due to the fact that $\RP(T)$ provides a useful ``alternative group of units'' in the case that the regular semigroup~$T$ is not a monoid.


Recall that there is a natural partial order $\preceq$ on the regular semigroup $T$: for $x,y\in T$, we say that $x\preceq y$ if and only if $x=ey=yf$ for some idempotents $e,f\in E(T)$; see for example \cite{Mitsch1986}.  It is easy to check that if $e,f\in E(T)$, then $e\preceq f$ if and only if $e=fef$.  We write $\MaxE(T)$ for the set of all $\preceq$-maximal idempotents of $T$.  It is easy to see that $\MI(T)\sub\MaxE(T)$.

Recall from \cite{BH1984} that the regular semigroup $T$ is \emph{RP-dominated} if every element of $T$ is $\preceq$-below an element of $\RP(T)$.  Analogously, we say that $T$ is \emph{MI-dominated} if every \emph{idempotent} of $T$ is $\preceq$-below an element of~$\MI(T)$.  (It is easy to see that any element $\preceq$-below a mid-identity must be an idempotent.)  It immediately follows that $\MaxE(T)=\MI(T)$ in any MI-dominated semigroup.

It will be convenient to state two technical results that we will need on a number of occations.  Part (i) of the next lemma is \cite[Lemma 2.5(1)]{BH1984}.  Part (ii) is \cite[Theorem 1.2]{BH1984}; see also \cite[Corollary 4.8]{Hickey1983}.

\begin{lemma}\label{lem:BH}
Let $T$ be a regular semigroup.
\bit
\itemit{i} If $x\in T$ and $e,f\in E(T)$ are such that $e\preceq x$ and $x\H f$, then $e\preceq f$.
\itemit{ii} If $T$ has a mid-identity, then $\RP(T)$ is a rectangular group and 
consists of those elements of $T$ that are $\H$-related to a mid-identity. \epfres
\end{itemize}
\end{lemma}

The next result 
shows that MI-domination is a weaker condition than RP-domination for regular semigroups with mid-identities.

\begin{prop}\label{prop:RP=>MI}
Suppose $T$ is a regular semigroup with a mid-identity.  If $T$ is RP-dominated, then $T$ is MI-dominated.
\end{prop}

\pf Suppose $T$ is RP-dominated, and let $e\in E(T)$.  We must show that $e\preceq u$ for some $u\in\MI(T)$.  By assumption, $e\preceq x$ for some $x\in\RP(T)$.  By Lemma \ref{lem:BH}(ii), $x\H u$ for some $u\in\MI(T)$.  By Lemma \ref{lem:BH}(i), since $e,u\in E(T)$, $e\preceq x$ and $x\H u$ together imply $e\preceq u$, as required. \epf

The converse of Proposition \ref{prop:RP=>MI} does not hold in general; for example, if $T$ is a regular monoid with trivial group of units and at least one non-idempotent.  We will soon give necessary and sufficient conditions for an MI-dominated regular semigroup to be RP-dominated; see Proposition \ref{prop:RP=>RP}.

Recall that if $e$ is an idempotent of a semigroup $T$, then the subsemigroup $eTe=\set{ete}{t\in T}$ of $T$ is a monoid with identity~$e$: the so-called \emph{local monoid} of $T$ with respect to $e$; these have already appeared in our investigations (see Remark \ref{rem:baSa}).  
It was shown in \cite[Lemma 3.2]{BH1984} that any RP-dominated regular semigroup with a mid-identity is the union of its local monoids corresponding to mid-identities.  Among other things, the next result generalises this to MI-dominated regular semigroups.

\begin{prop}\label{prop:MI_LMC}
Let $T$ be a regular semigroup, write $R=\RP(T)$ and $M=\MI(T)$, and suppose $M\not=\emptyset$.
\bit
\itemit{i} If $e\in M$, then the map $T\to eTe:x\mt exe$ is an epimorphism.
\itemit{ii} If $e,f\in M$, then the maps $eTe\to fTf:x\mt fxf$ and $fTf\to eTe:x\mt exe$ are mutually inverse isomorphisms.
\itemit{iii} The set $\bigcup_{e\in M}eTe = MTM = RTR$ is a subsemigroup of $T$.
\itemit{iv} $T$ is MI-dominated if and only if $T=\bigcup_{e\in M}eTe$.
\end{itemize}
\end{prop}

\pf (i) and (ii).  These are easily checked.

\pfitem{iii}  It is clear that $RTR$ is a subsemigroup of $T$.  It is also clear that $\bigcup_{e\in M}eTe \sub MTM \sub RTR$.  Now suppose $x,z\in R$ and $y\in T$; we must show that $xyz\in eTe$ for some $e\in M$.  By Lemma \ref{lem:BH}(ii), $x\H u$ and $z\H v$ for some $u,v\in\MI(T)$.  In particular, $x=ux$ and $z=zv$.  Since $u,v\in\MI(T)$, it follows that $xyz=(ux)y(zv)=(uv)xyz(uv)\in eTe$, where $e=uv\in\MI(T)$. 

\pfitem{iv}  Suppose $T$ is MI-dominated, and let $x\in T$.  Since $T$ is regular, $x=exf$ for some $e,f\in E(T)$.  By assumption, $e\preceq u$ and $f\preceq v$, for some $u,v\in\MI(T)$.  In particular, $e=ue$ and $f=fv$, and so $x=exf=u(exf)v\in MTM=\bigcup_{e\in M}eTe$, using (iii) in the last step.  

Conversely, suppose $T=\bigcup_{e\in M}eTe$, and let $u\in E(T)$.  Then $u\in eTe$ for some $e\in M$.  It immediately follows that $u\preceq e$. \epf

So MI-dominated regular semigroups are precisely those regular semigroups that are ``covered'' (in some non-technical sense) by local monoids corresponding to mid-identities. 
Before we give the criteria for an MI-dominated regular semigroup to be RP-dominated, we first consider the situation of RP-dominated \emph{monoids}.

Suppose $T$ is a monoid with identity $e$.  Clearly $e\in\MI(T)$.  On the other hand, if $u\in\MI(T)$, then $e=ee=eue=u$, so it follows that $\MI(T)=\{e\}$.  Clearly all idempotents of $T$ are $\preceq$-below $e$, and so $T$ is MI-dominated.  By Lemma \ref{lem:BH}(ii), $\RP(T)=H_e$ is the $\H$-class of $e$: that is, $\RP(T)$ is the group of units of $T$ (this was also shown in \cite[Proposition 1]{KL2001}).  Thus, the monoid $T$ is RP-dominated if and only if every element of $T$ is $\preceq$-below a unit of $T$.  It is also easy to see that if $x\in T$ and $y\in H_e$, then
\[
x\preceq y \iff x=fy\ (\exists f\in E(T)) \iff x=yf\ (\exists f\in E(T)).
\]
Thus, a monoid is RP-dominated if and only if it is \emph{factorisable} in the sense of \cite{Tirasupa1979,CH1974}.

\begin{prop}\label{prop:RP=>RP}
Let $T$ be an MI-dominated regular semigroup.  Then $T$ is RP-dominated if and only if the local monoid $eTe$ is RP-dominated (equivalently, factorisable) for each mid-identity $e\in\MI(T)$.
\end{prop}

\pf Suppose first that $T$ is RP-dominated, and let $e\in\MI(T)$.  It is easy to show that $eTe$ is regular.  In order to show that $eTe$ is RP-dominated, let $x\in eTe$.  Since $T$ is RP-dominated, we have $x\preceq y$ (in~$T$) for some $y\in\RP(T)$.  So $x=fy=yg$ for some $f,g\in E(T)$.  Since $e\in\MI(T)$ and $x\in eTe$, we have $x=exe=e(fy)e=(efe)(eye)$.  A similar calculation gives $x=(eye)(ege)$.  It is easy to check that $efe,ege\in E(eTe)$ and $eye\in\RP(eTe)$, so it follows that $x\preceq eye$ (in $eTe$).  This shows that $eTe$ is RP-dominated.

Conversely, suppose $eTe$ is RP-dominated for each $e\in\MI(T)$.  Let $x\in T$ be arbitrary.  Since $T$ is MI-dominated, Proposition \ref{prop:MI_LMC}(iv) gives $x\in eTe$ for some $e\in\MI(T)$.  Since $eTe$ is RP-dominated, we have $x\preceq y$ (in $eTe$, and hence also in $T$) for some $y\in\RP(eTe)$.  But $eTe$ is a monoid, and so $\RP(eTe)$ is the group of units of $eTe$; thus, $y$ is $\H$-related (in $eTe$, and hence also in $T$) to $e$.  It follows from Lemma \ref{lem:BH}(ii) that $y\in\RP(T)$, and this competes the proof that $T$ is RP-dominated. \epf

\subsectiontitle{Mid-identities and regularity-preserving elements in $P^a$}

We now bring our focus back to the regular semigroup $P^a=\RegSija$.  Proposition \ref{prop:MI_LMC} clearly applies to $T=P^a$, since $b\in\MI(P^a)$.
The next result describes the mid-identities and regularity-preserving elements in $P^a$.  Recall that $W=(aS_{ij}a,\starb)$ is a regular monoid with identity $a$, where $\starb$ is the restriction of the $\star_b$ operation to $aS_{ij}a\sub S_{ji}$.


\begin{prop}\label{prop:MI_RP}
We have 
\begin{itemize}\begin{multicols}{2}
\itemit{i} $\MI(P^a)=E_a(\Hh_b^a)=V(a)\sub\MaxE(P^a)$, 
\itemit{ii} $\RP(P^a)=\Hh_b^a$.  
\end{multicols}\end{itemize}
\end{prop}

\pf (i).  Suppose $u\in\MI(P^a)$.  Since $\MI(P^a)\sub E_a(P^a)$, as noted above, we have $u=uau$.  Since also
\[
a=ababa=a(b\star_ab)a=a(b\star_au\star_ab)a=abauaba=aua,
\]
we have $u\in V(a)$.  This shows that $\MI(P^a)\sub V(a)$.  

Next, suppose $u\in V(a)$.  So $u=uau$ and $a=aua$.  The former shows that $u$ is an idempotent of $P^a$.  Using the latter, we see that for any $x\in P^a$, $\ub\starb\xb=auaxa=axa=\xb$, and similarly $\xb\starb\ub=\xb$, so that $\ub=a=\bb$ is the identity of $W$.  Thus, $u\gHh^ab$, so that $u\in\Hh_b^a$.  Since we have already observed that $u$ is an idempotent, it follows that $u\in E_a(\Hh_b^a)$.  Thus, $V(a)\sub E_a(\Hh_b^a)$.

Finally, suppose $u\in E_a(\Hh_b^a)$.  So $\ub=\bb=a$ is the unique idempotent of $H_{\bb}^b=H_a^b$.  That is, $aua=a$.  It follows that for any $x,y\in P^a$, $x\star_au\star_ay=x(aua)y=xay=x\star_ay$, giving $u\in\MI(P^a)$.  So $E_a(\Hh_b^a)\sub\MI(P^a)$.

We have already noted that $\MI(T)\sub\MaxE(T)$ in any regular semigroup $T$.

\pfitem{ii}  Since $P^a$ is regular and $\MI(P^a)\not=\emptyset$, Lemma \ref{lem:BH}(ii) says that $\RP(P^a)=\bigcup_{e\in\MI(P^a)}H_e^a$.  Since $\MI(P^a)=E_a(\Hh_b^a)$, by part (i), and since every $\H^a$-class contained in $\Hh_b^a$ is a group, it follows that $\RP(P^a)=\bigcup_{e\in E_a(\Hh_b^a)} H_e^a=\Hh_b^a$. \epf

%

Before we move on, it will be convenient to prove the following result, which gives an alternative description of $\Hh_b^a=\RP(P^a)$ in a certain special case.

\begin{lemma}\label{lem:Jab1}
If $J_a^\starb =H_a^\starb $ in $W$, then $\Hh_b^a=D_b^a=J_b^{P^a}$ in $P^a$.
\end{lemma}

\pf Since ${\gHh}^a\sub{\D^a}\sub{\JPa}$, by Lemma \ref{lem:K_hat}, we have $\Hh_b^a\sub D_b^a\sub J_b^{P^a}$.  Let $x\in J_b^{P^a}$; we must show that $x\in\Hh_b^a$.  Now, $x\in J_b^{P^a}$ gives $\xb\in J_{\bb}^\starb =J_a^\starb =H_a^\starb $.  So $\xb\H^\starb a=\bb$, giving $x\gHh^ab$, and $x\in\Hh_b^a$. \epf

In particular, if $W$ is stable (which occurs if $W$ is finite, for example, by \cite[Theorem A.2.4]{RSbook}), then $J_a^\starb =H_a^\starb $, as the next lemma shows.

%

\begin{lemma}\label{lem:Jab2}
Let $M$ be a monoid with identity $e$.
\bit
\itemit{i} If every element of $J_e$ is $\R$-stable, then $J_e=R_e$.
\itemit{ii} If every element of $J_e$ is $\L$-stable, then $J_e=L_e$.
\itemit{iii} If every element of $J_e$ is stable, then $J_e=H_e$.
\end{itemize}
\end{lemma}

\pf We just prove (i), since (ii) is dual, and (iii) follows from (i) and (ii), as $H_e=R_e\cap L_e$.  Suppose every element of $J_e$ is $\R$-stable.  Clearly $R_e\sub J_e$.  Conversely, if $x\in J_e$, then $e\J x=ex$; $\R$-stability of $x$ then gives $e\R ex=x$, so $x\in R_e$.  \epf

For $e\in\MI(P^a)=V(a)$, we write $W_e$ for the local monoid $e\star_aP^a\star_ae$ of $P^a$.  
%
By parts (i) and (ii) of Proposition \ref{prop:MI_LMC}, we know that all of these local monoids are isomorphic to each other, and are homomorphic images of $P^a$ itself; the next result takes this further, and also explains our choice of notation.

\begin{prop}\label{prop:We_new}
For any $e\in V(a)$, the restriction of $\phi$ to $W_e$ is an isomorphism $\phi|_{W_e}:W_e\to W$.
\end{prop}

\pf It is easy to check that the map $W\to W_e:x\mt exe$ is the inverse of $\phi|_{W_e}$. \epf


Thus, if $P^a=\RegSija$ is MI-dominated, then Propositions \ref{prop:MI_LMC}(iv) and \ref{prop:We_new} say that $P^a$ is a union of isomorphic copies of $W=(aS_{ij}a,\starb)$.  All of the examples we consider in \cite{Sandwiches2} satisfy the MI-domination property; for some natural examples that do not, see \cite{DEdiagrams}.

\subsectiontitle{Rank and idempotent rank}

Recall that the \emph{rank} of a semigroup $T$ is the cardinal
\begin{align*}
\rank(T)&=\min\bigset{|A|}{A\sub T,\ \la A\ra=T}.
\intertext{If~$T$ is idempotent-generated, then the \emph{idempotent rank} of $T$ is}
\idrank(T)&=\min\bigset{|A|}{A\sub E(T),\ \la A\ra=T}.
\end{align*}
As noted above, we wish to obtain formulae for the rank of $P^a$, and for the rank and idempotent rank of~$\E_a(P^a)$.
Before this, we must first prove some preliminary results.  If $M$ is a monoid, we will write $G_M$ for the \emph{group of units} of $M$: that is, the $\H$-class of the identity element of $M$.  
%
Recall from \cite{HHR} that the \emph{relative rank} of a semigroup $T$ with respect to a subset $A\sub T$ is
\begin{align*}
\rank(T:A)&=\min\bigset{|B|}{B\sub T,\ \la A\cup B\ra=T}.
\intertext{If $T$ is idempotent-generated, then the \emph{relative idempotent rank} of $T$ with respect to a subset $A\sub E(T)$ is}
\idrank(T:A)&=\min\bigset{|B|}{B\sub E(T),\ \la A\cup B\ra=T}.
\end{align*}

\begin{lemma}\label{lem:GMid}
Suppose $M$ is an idempotent-generated monoid.  Then 
\begin{itemize}\begin{multicols}{2}
\itemit{i} $G_M=\{\id_M\}$,
\itemit{ii} $M\sm G_M$ is an ideal of $M$,
\itemit{iii} $\rank(M)=1+\rank(M:G_M)$,
\itemit{iv} $\idrank(M)=1+\idrank(M:G_M)$.
\end{multicols}\eitmc
\end{lemma}

\pf (i).  Clearly $\{\id_M\}\sub G_M$.  Conversely, suppose $g\in G_M$, and consider an expression $g=e_1\cdots e_k$, where ${e_1,\ldots,e_k\in E(M)}$ and $k$ is minimal.  Then $g=e_1g$, so $\id_M=gg^{-1}=e_1gg^{-1}=e_1\id_M=e_1$.  If $k\geq2$, then we would have $g=\id_M\cdot e_2\cdots e_k=e_2\cdots e_k$, contradicting the minimality of $k$.  So $k=1$, and $g=e_1=\id_M$, as required.


\pfitem{ii}  Suppose $x\in M\sm G_M$ and $y\in M$.  We must show that $xy,yx\in M\sm G_M$.  We just do this for $xy$, as the proof for $yx$ is dual.  Suppose to the contrary that $xy\in G_M=\{\id_M\}$.  Write $x=e_1\cdots e_k$, where $e_1,\ldots,e_k\in E(M)$ and $k$ is minimal.  Then $x=e_1x$, so that $\id_M=xy=e_1xy=e_1\id_M=e_1$.  
As in the previous paragraph, we deduce that $k=1$, in which case $x=e_1=\id_M$, contradicting $x\in M\sm G_M$.


\pfitem{iii) and (iv}  These follow immediately from 
the fact that any generating set for $M$ must contain $\id_M$, which itself follows from the fact that $M\sm G_M$ is an ideal.  \epf

The next technical result is pivotal in what follows.  It may be regarded as a strengthening of Lemma~\ref{lem:Xb} under a certain condition that is weaker than MI-domination.

\begin{lemma}\label{lem:XY}
Suppose every idempotent of $P^a$ is $\preceq$-below a maximal one, and suppose $X\sub P^a$ is such that $E_b(W)\sub\la\Xb\ra_b$ and $\MaxE(P^a)\sub\la X\ra_a$.  Then $\la\Xb\ra_b\phi^{-1}=\la X\ra_a$.
\end{lemma}

\pf Clearly $\la X\ra_a\sub\la\Xb\ra_b\phi^{-1}$.  Conversely, by Lemma \ref{lem:Xb}, we know that $\la\Xb\ra_b\phi^{-1}\sub\la X\cup E_a(P^a)\ra_a$, so it suffices to show that $E_a(P^a)\sub \la X\ra_a$.  With this in mind, let $e\in E_a(P^a)$.  Then $e\preceq f$ for some $f\in\MaxE(P^a)$.  So $e=f\star_ae\star_af=faeaf$.  Now, $\eb\in E_b(W)\sub\la\Xb\ra_b$, so we may write $\eb=\xb_1\starb\cdots\starb\xb_k$ for some $x_1,\ldots,x_k\in X$: that is, $aea=ax_1ax_2a\cdots ax_ka$.  But then $e=faeaf=fax_1ax_2a\cdots ax_kaf=f\star_ax_1\star_ax_2\star_a\cdots\star_ax_k\star_af$.  Since $f\in\MaxE(P^a)\sub\la X\ra_a$, it follows that $e\in\la X\ra_a$, as required. \epf


We will use the following result of Ru\v skuc \cite{Ruskuc1994} several times in our subsequent investigations; the various statements are spread across several results and proofs (concerning the more general \emph{completely 0-simple semigroups}) from \cite{Ruskuc1994}.  Thus, for convenience, we give a short proof in the special case of rectangular~groups.

\begin{prop}\label{prop:Ruskuc}
Let $T$ be an $r\times l$ rectangular group over $G$.  Then
\bit
\itemit{i} $\rank(T)=\max(r,l,\rank(G))$,
\itemit{ii} any generating set for $T$ contains elements from every $\R$-class, and from every $\L$-class, of $T$,
\itemit{iii} if $\rank(T)=r$, then there is a minimum-size generating set for $T$ that is a cross-section of the $\R$-classes of~$T$,
\itemit{iv} if $\rank(T)=l$, then there is a minimum-size generating set for $T$ that is a cross-section of the $\L$-classes of~$T$.
\eitres
\end{prop}

\pf Throughout the proof, we assume $T=I\times G\times J$, with multiplication $(i_1,g,j_1)(i_2,h,j_2)=(i_1,gh,j_2)$, and with $r=|I|$ and $l=|J|$.  We fix a generating set $\Ga$ of $G$ with $|\Ga|=\rank(G)$.  We also fix some set $X$ with $|X|=\max(r,l,\rank(G))$, and surjective mappings
\[
\al:X\to I:x\mt i_x \COMMA
\be:X\to J:x\mt j_x \COMMA
\ga:X\to \Ga:x\mt g_x ,
\]
assuming that $\al$ (or $\be$ or $\ga$) is a bijection if $r=|X|$ (or $l=|X|$ or $\rank(G)=|X|$, respectively).  For $x\in X$, put $t_x=(i_x,g_x,j_x)$.  We first claim that $T=\la\Om\ra$, where $\Om=\set{t_x}{x\in X}$.  Indeed, suppose $i\in I$, $j\in J$ and $g\in G$ are arbitrary.  Then we may write $i=i_x$, $j=j_y$ and $g_x^{-1}gg_y^{-1}=g_{z_1}\cdots g_{z_k}$, for some $x,y,z_1,\ldots,z_k\in X$.  But then $(i,g,j)=t_xt_{z_1}\cdots t_{z_k}t_y$, establishing the claim.  

\pfitem{ii}  Suppose $T=\la \Upsilon\ra$, and let $t=(i,g,j)\in T$ be arbitrary.  Then $R_t=\{i\}\times G\times J$ and $L_t=I\times G\times\{j\}$.  Now consider an expression $t=(i_1,g_1,j_1)\cdots(i_k,g_k,j_k)$, where the factors belong to $\Upsilon$.  Then $t=(i_1,g_1\cdots g_k,j_k)$, so that $i=i_1$ and $j=j_k$, giving $(i_1,g_1,j_1)\in\Upsilon\cap R_t$ and $(i_k,g_k,j_k)\in\Upsilon\cap L_t$.

\pfitem{i}  By the first paragraph, $\rank(T)\leq|\Om|=|X|=\max(r,l,\rank(G))$.  By (ii), $\rank(T)\geq|T/{\R}|=|I|$, and similarly $\rank(T)\geq|J|$.  Since the map $T\to G:(i,g,j)\mt g$ is an epimorphism, we also have $\rank(T)\geq\rank(G)$.

\pfitem{iii) and (iv}  These both follow from the first paragraph.  For example, if $r=\rank(T)=|X|$, then the fact that $\al$ is a bijection shows that $\Om$ is a cross-section of $T/{\R}$. \epf


The next general result will allow us to quickly derive the promised results concerning the rank of $P^a$ and the (idempotent) rank of $\E_a(P^a)$.
Recall that a subsemigroup of a semigroup $T$ is \emph{full} if it contains $E(T)$.  It is easy to show that any full subsemigroup of a rectangular group $I\times G\times J$ is of the form $I\times K\times J$, where $K$ is some submonoid of $G$.

%
%

\begin{prop}\label{prop:MN}
Let $r=|\Hh_b^a/{\R^a}|$ and $l=|\Hh_b^a/{\L^a}|$, let $M$ be a full submonoid of~$W$ for which $M\sm G_M$ is an ideal of $M$, and $G_M=M\cap G_W$, and put $N=M\phi^{-1}$.
Then 
\[
\rank(N)\geq\rank(M:G_M)+\max(r,l,\rank(G_M)),
\]
with equality if $P^a$ is MI-dominated.
\end{prop}

\pf First, note that $N$ is a full subsemigroup of $P^a$, by Lemma \ref{lem:EaEb}.  In particular, $N\cap\Hh_b^a$ is a full subsemigroup of $\Hh_b^a$, and we recall that $\Hh_b^a$ is an $r\times l$ rectangular group over $H_b^a\cong H_a^\starb=G_W$, by Theorem~\ref{thm:RG}(iii).  
From the assumption that $G_M=M\cap G_W$, and the above observation about full subsemigroups of rectangular groups, it quickly follows that $N\cap\Hh_b^a$ is an $r\times l$ rectangular group over some subgroup $K$ of $H_b^a$.  In other words, $N\cap\Hh_b^a=V(a)\star_a K\star_aV(a)$; recall that $V(a)=E_a(\Hh_b^a)$.  But $K\phi=(N\cap\Hh_b^a)\phi=M\cap H_a^b=G_M$.  
Since $\phi|_{H_b^a}$ is injective (by Theorem~\ref{thm:RG}(ii)), so too is $\phi|_K$, and it follows that $K\cong G_M$.
Since $N\cap\Hh_b^a$ is a rectangular group, it then follows from Proposition \ref{prop:Ruskuc} that
\[
\rank(N\cap\Hh_b^a)=\max(r,l,\rank(K))=\max(r,l,\rank(G_M)).
\]
Now suppose $N=\la X\ra_a$, with $|X|=\rank(N)$, and put $Y=X\cap\Hh_b^a$ and $Z=X\sm \Hh_b^a$.  Let $u\in N\cap\Hh_b^a$, and consider an expression $u=x_1\star_a\cdots\star_ax_k$, where $x_1,\ldots,x_k\in X$.  Then $\ub=\xb_1\starb\cdots\starb\xb_k$.  Since $\ub\in G_M$, and since $M\sm G_M$ is an ideal of $M$, it follows that $\xb_1,\ldots,\xb_k\in G_M$.  Consequently, $x_1,\ldots,x_k\in G_M\phi^{-1}\sub\Hh_b^a$, giving $x_1,\ldots,x_k\in Y$.  It follows that $N\cap\Hh_b^a=\la Y\ra_a$.  Consequently, 
\begin{equation}\label{eq:Y}
|Y|\geq\rank(N\cap\Hh_b^a)
=\max(r,l,\rank(G_M)).
\end{equation}
Next, we note that $M=N\phi=\la X\ra_a\phi=\la\Xb\ra_b=\la\Yb\cup\Zb\ra_b=\la\la\Yb\ra_b\cup\Zb\ra_b=\la G_M\cup\Zb\ra_b$, so that 
\begin{equation}\label{eq:Z}
|Z|\geq|\Zb|\geq\rank(M:G_M).
\end{equation}
Equations \eqref{eq:Y} and \eqref{eq:Z} then give $\rank(N)=|X|=|Y|+|Z|\geq\max(r,l,\rank(G_M))+\rank(M:G_M)$, as required.

For the remainder of the proof, suppose $P^a$ is MI-dominated.  It suffices to give a generating set of the stated size.  With this in mind, suppose $N\cap\Hh_b^a=\la Y\ra_a$ with
$|Y|=\rank(N\cap\Hh_b^a)=\max(r,l,\rank(G_M))$.  Also, let $Z\sub P^a$ be such that $M=\la G_M\cup\Zb\ra_b$ and $|Z|=\rank(M:G_M)$.
Put $X=Y\cup Z$.  Since $X$ has the desired size,
by construction, it suffices to show that $N=\la X\ra_a$.  Now,
\[
M=\la G_M\cup\Zb\ra_b = \la(N\cap\Hh_b^a)\phi\cup\Zb\ra_b = \big\la \la Y\ra_a\phi\cup\Zb\big\ra_b = \big\la\la\Yb\ra_b\cup\Zb\big\ra_b = \la\Yb\cup\Zb\ra_b=\la\Xb\ra_b.
\]
Also note that $E_b(W)\sub M=\la\Xb\ra_b$, as $M$ is full, and that
\[
\MaxE(P^a)=\MI(P^a)=V(a)=E_a(\Hh_b^a)\sub N\cap\Hh_b^a=\la Y\ra_a\sub\la X\ra_a.
\]
It follows from Lemma \ref{lem:XY} that $N=M\phi^{-1}=\la\Xb\ra_b\phi^{-1}=\la X\ra_a$.  As noted above, this completes the proof.~\epf

%
%
%

\begin{rem}
In the statement of Proposition \ref{prop:MN} (see also Theorems \ref{thm:rankU} and \ref{thm:rankE}), the condition ``$P^a$ is MI-dominated'' is equivalent to ``$N$ is MI-dominated''.  Indeed, this follows from the following claim (and the fact that idempotents satisfy $e\preceq f \iff e=fef$ in any (regular) semigroup), since~$N$ is full (as observed at the beginning of the above proof):
\bit
\item[] \emph{If $U$ is a full subsemigroup of a regular semigroup $T$ with a mid-identity, then $\MI(U)=\MI(T)$.}
\eit
To establish the claim, first suppose $u\in\MI(T)$.  In particular, $u\in E(T)\sub U$, and it immediately follows that $u\in\MI(U)$.  Conversely, suppose $u\in\MI(U)$, and let $e\in\MI(T)$ be arbitrary.  Since $e\in E(T)\sub U$, it follows that $e=ee=eue$.  Then for any $x,y\in T$, $xy=xey=xeuey=xuy$, giving $u\in\MI(T)$.
\end{rem}

The two most obvious full submonoids of $W$ are $W$ itself, and $\E_b(W)$, the idempotent-generated submonoid of $W$.  Recall that $W$ is a monoid with group of units $G_W=H_a^b$.  The next result follows immediately from Proposition \ref{prop:MN}, with $M=W$, noting that $W\phi^{-1}=P^a$, as $\phi$ is surjective.

\begin{thm}\label{thm:rankU}
Let $r=|\Hh_b^a/{\R^a}|$ and $l=|\Hh_b^a/{\L^a}|$, and suppose $W\sm G_W$ is an ideal of $W$.  Then
\[
\rank(P^a)\geq\rank(W:G_W)+\max(r,l,\rank(G_W)),
\]
with equality if $P^a$ is MI-dominated. \epfres
\end{thm}

The statement for the corresponding result on $\E_a(P^a)=\E_a(\Sija)$ is somewhat simpler, since the group of units of $\E_b(W)$ is $\{a\}$, and $\E_b(W)\sm\{a\}$ is an ideal,
by Lemma \ref{lem:GMid}.

%
%

\begin{thm}\label{thm:rankE}
Let $r=|\Hh_b^a/{\R^a}|$ and $l=|\Hh_b^a/{\L^a}|$.  Then
\[
\rank(\E_a(P^a)) \geq \rank(\E_b(W))+\max(r,l)-1 \AND \idrank(\E_a(P^a)) \geq \idrank(\E_b(W))+\max(r,l)-1,
\]
with equality in both if $P^a$ is MI-dominated.
\end{thm}

\pf Put $M=\E_b(W)$, so that $N=M\phi^{-1}=\E_a(P^a)$, by Theorem \ref{thm:EaEb}.  By Lemma \ref{lem:GMid}, $M\sm G_M$ is an ideal of $M$, so the conditions of Proposition \ref{prop:MN} are satisfied, and we obtain 
\[
\rank(\E_a(P^a)) =\rank(N) \geq \rank(M:G_M)+\max(r,l,\rank(G_M)),
\]
with equality in the MI-dominated case.
Lemma \ref{lem:GMid} gives $\rank(M:G_M)=\rank(M)-1$, and also ${G_M=\{\id_M\}=\{a\}}$; from the latter, we immediately obtain $\rank(G_M)=1$, and so $\max(r,l,\rank(G_M))=\max(r,l)$.

Next we prove the assertion about $\idrank(N)$.  Since $G_M=\{a\}$, we have $N\cap\Hh_b^a=a\phi^{-1}=V(a)=E_a(\Hh_b^a)$.  Suppose $M=\la X\ra_a$, where $X\sub E_a(P^a)$, and put $Y=X\cap V(a)$ and $Z=X\sm V(a)$.  As in the proof of Proposition \ref{prop:MN}:
\bit
\item[(i)] $V(a)=\la Y\ra_a$, so $|Y|\geq\rank(V(a))=\max(r,l)$,
\item[(ii)] $M=\la G_M\cup\Zb\ra_b=\la\{a\}\cup\Zb\ra_b$, so that $|Z|\geq|\Zb|\geq\idrank(M:G_M)=\idrank(M)-1$, 
\item[(iii)] if $P^a$ is MI-dominated, then $N=\la Y_1\cup Z_1\ra_a$ if $Y_1,Z_1\sub E_a(P^a)$ are such that $V(a)=\la Y_1\ra_a$ and $\la\{a\}\cup\Zb_1\ra_b=M$, with $|Y_1|=\max(r,l)$ and $|Z_1|=\idrank(M:G_M)$.
\eit
It follows from (i) and (ii) that $\idrank(N)\geq\idrank(\E_b(W))-1+\max(r,l)$, while (iii) leads to the converse in the MI-dominated case.~\epf

\begin{rem}
While MI-domination is sufficient to obtain equality in Theorems \ref{thm:rankU} and \ref{thm:rankE}, it is not \emph{necessary}; see \cite{DEdiagrams}.
\end{rem}

\begin{rem}\label{rem:rankE}
If $P^a$ is MI-dominated, then it follows from Theorem \ref{thm:rankE} that
\[
\rank(\E_b(W))=\idrank(\E_b(W)) \implies \rank(\E_a(P^a))=\idrank(\E_a(P^a)).
\]
The converse holds if $r,l<\aleph_0$.
\end{rem}

%

\sectiontitle{Inverse monoids}\label{sect:inverse}

The general theory of $P^a=\RegSija$ developed above in Sections \ref{sect:RegSija} and \ref{sect:LMC} relied on the assumption that the sandwich element $a\in S_{ji}$ is \emph{sandwich-regular}: that is, every element of $aS_{ij}a$ is regular (in $S$).  
In this short section, we prove a result that shows how the general theory simplifies substantially in the case that $a$ satisfies a certain stronger version of sandwich-regularity.

Recall that a semigroup $T$ is \emph{inverse} if, for all $x\in T$, there exists a unique $y\in T$ such that $x=xyx$ and $y=yxy$.  
In particular, inverse semigroups are regular.
If $x$ is an element of the partial semigroup $S\equiv(S,\cdot,I,\lam,\rho)$, we say that $x$ is \emph{uniquely regular} if the set $V(x)=\set{y\in S}{x=xyx,\ y=yxy}$ has size~$1$.  We say that $a\in S_{ji}$ is \emph{uniquely sandwich-regular} if every element of $\{a\}\cup aS_{ij}a$ is uniquely regular (in~$S$).  
For example, every element is uniquely regular (and, hence, uniquely sandwich-regular) in a (locally small) \emph{inverse category}, as defined in \cite{Kastl1979}; see also \cite[Section 2.3.2]{CL2002}.

\begin{prop}\label{prop:inverse}
Suppose $a\in S_{ji}$ is uniquely sandwich-regular, and that $V(a)=\{b\}$.  Then all maps in the diagram \eqref{eq:CD_Reg} are isomorphisms, and all semigroups are inverse monoids.
\end{prop}

\pf We first show that $P^a=\RegSija$ is a monoid with identity $b$.  Let $x\in P^a$, and let $y\in V_a(x)$.
It is routine to check that $x,xab,bax\in V(aya)$.  
Since $aya\in aS_{ij}a$ is uniquely regular, it follows that $x=xab=bax$: that is, $x=x\star_ab=b\star_ax$, as required.


Next we show that $\phi$ is injective.  Since $\phi=\psi_1\phi_1=\psi_2\phi_2$ and all maps in \eqref{eq:CD_Reg} are epimorphisms, it will then follow that all maps are injective, and hence isomorphisms.  So suppose $x,y\in P^a$ are such that $axa=x\phi=y\phi=aya$.  Then, since $b$ is the identity element of $P^a$, $x=baxab=bayab=y$.

To complete the proof, it suffices to show that $P^a$ is inverse.  We know that $P^a$ is regular.  Suppose now that $u\in P^a$ and that $x,y\in V_a(u)$.  It follows quickly that $x,y\in V(aua)$ so, since $aua\in aS_{ij}a$ is uniquely regular, it follows that $x=y$, as required. \epf

\begin{rem}\label{rem:inverse}
Note that many of the results of Sections \ref{sect:RegSija} and \ref{sect:LMC} become trivial, or even vacuous, if $a\in S_{ji}$ is uniquely sandwich-regular.  For example:
\bit
\item the map $\psi=(\psi_1,\psi_2)$ from Theorem \ref{thm:psi} is trivially injective,
\item the $\gKh^a$ relations of Section \ref{sect:Gh} are identical to the $\K^a$ relations, so the rectangular groups in Theorem~\ref{thm:RG} are just groups,
\item Theorem \ref{thm:EaEb} is completely trivial,
\item Proposition \ref{prop:MI_RP} says that $\MI(P^a)=\{b\}$ and $\RP(P^a)=H_b^a$ consist only of the identity and the units, respectively, which is true in any monoid,
\item Theorem \ref{thm:rankU} reduces to ``$\rank(W)=\rank(W:G_W)+\rank(G_W)$ if $W\sm G_W$ is an ideal of $W$'',
\item Theorem \ref{thm:rankE} reduces to ``$\rank(\E_a(P^a))=\rank(\E_b(W))$ and $\idrank(\E_a(P^a))=\idrank(\E_b(W))$''.
\eit
\end{rem}

\sectiontitle{The rank of a sandwich semigroup}\label{sect:rankSija}

We conclude with another short section; this one contains some remarks about the rank
of a sandwich semigroup $\Sija$.  Unlike the case with the regular and idempotent-generated subsemigroups, $\RegSija$ and~$\E_a(\Sija)$, we are not able to give precise values for $\rank(\Sija)$, even under assumptions such as sandwich-regularity of~$a$, though we can state some rough lower bounds that will be useful in our concrete applications in \cite{Sandwiches2}.
%
We fix a partial semigroup $S\equiv(S,\cdot,I,\lam,\rho)$, and a sandwich element $a\in S_{ji}$ for some $i,j\in I$, but we make no further assumptions (such as regularity of $a$).  

Recall that we have orders $\leq_{\J}$ and $\leq_{\J^a}$ on the sets $S_{ij}/{\J}$ and $S_{ij}/{\J^a}$ of all $\J$- and $\J^a$-classes in~$S_{ij}$.  It is clear that any generating set for $\Sija$ must include elements from every maximal $\J^a$-class (consider a factorisation over the generators of an element from such a maximal $\J^a$-class).  
The next result concerns maximal $\J$-classes, and provides more precise information under certain stability assumptions that will occur in our concrete examples.

\begin{lemma}\label{lem:rankSija2}
Suppose $\Sija=\la X\ra_a$, and that $J$ is a maximal $\J$-class in $S_{ij}$.  Then 
\bit
\itemit{i} $X\cap J\not=\emptyset$,
\itemit{ii} if every element of $aS_{ij}$ is $\R$-stable, then $X$ has non-trivial intersection with each $\R$-class of $J$, 
\itemit{iii} if every element of $S_{ij}a$ is $\L$-stable, then $X$ has non-trivial intersection with each $\L$-class of $J$.
\eitres
\end{lemma}

\pf (i).  Let $x\in J$, and consider an expression $x=x_1\star_a\cdots\star_ax_k$, where $x_1,\ldots,x_k\in X$.  Note that $x\leq_{\J}x_i$ for all $i$.  Since $J_x=J$ is maximal, it follows that $x\J x_i$ for all $i$.  In particular, each $x_i$ belongs to $X\cap J$.

\pfitem{ii}  Suppose every element of $aS_{ij}$ is $\R$-stable.  To prove (ii), it suffices to prove that $x_1\R x$.  This is clear if $k=1$, so suppose $k\geq2$, and write $z=x_2\star_a\cdots\star_ax_k$, so that $x=x_1\star_az$.  From above, $x_1\J x=x_1(az)$.  Since $az\in aS_{ij}$ is $\R$-stable, it follows that $x_1\R x_1(az)=x$, as required. Part (iii) is dual. \epf

In particular, if every element of $aS_{ij}$ is $\R$-stable, every element of $S_{ij}a$ is $\L$-stable, and the maximal $\J$-classes of $S_{ij}$ are $\set{J_k}{k\in K}$, then Lemma \ref{lem:rankSija2} gives
\[
\rank(\Sija)\geq\sum_{k\in K}\max\big(|J_k/{\R}|,|J_k/{\L}|\big).
\]
In several examples we consider in \cite{Sandwiches2}, there will be a single maximal $\J$-class in $S_{ij}$.  In some cases, but not all, the above lower bound is the exact value of $\rank(\Sija)$; see for example \cite[Theorems~2.63 and~2.69]{Sandwiches2}.

\section*{Acknowledgements}

This work was initiated during visits of the third author to Chiang Mai in 2015, and to Novi Sad in 2016; he thanks these institutions for their generous hospitality.  
The first and second authors are supported by Grant Nos.~174019 and 174018, respectively, of the Ministry of Education, Science, and Technological Development of the Republic of Serbia.
The first and third authors would like to thank \emph{Macchiato Fax} for providing us with many a \emph{Klub sendvi\v{c}}.
Finally, we thank the referee for his/her careful reading of a lengthy pair of articles, and for a number of helpful suggestions.

\footnotesize
\def\bibspacing{-1.1pt}
\bibliography{biblio}
\bibliographystyle{plain}
\end{document}